\documentclass[10pt,leqno]{article}
\usepackage{amsmath,amssymb,mathrsfs,graphicx,subfig,wrapfig}
\usepackage[affil-it]{authblk}
\newtheorem{thm}{Theorem}
\newtheorem{lemma}{Lemma}
\newtheorem{proposition}{Proposition}
\newtheorem{coro}{Corollary}

\def\XXint#1#2#3{{\setbox0=\hbox{$#1{#2#3}{\int}$ }
\vcenter{\hbox{$#2#3$ }}\kern-.6\wd0}}

\title{Multiplicity of solutions for polyharmonic Dirichlet problems with exponential nonlinearities and broken symmetry}

\author{Edger Sterjo \thanks{Electronic address: \texttt{esterjo@gradcenter.cuny.edu}}}

\affil{Department of Mathematics, The Graduate Center, CUNY, \\365 Fifth Avenue
New York, NY 10016 USA}

\date{}

\begin{document}

\maketitle

\begin{abstract}
We prove the existence of infinitely many solutions to a class of non-symmetric Dirichlet problems with exponential nonlinearities. Here the domain $\Omega \subset\subset \mathbb{R}^{2l}$ where $2l$ is also the order of the equation. Considered are the problem with no symmetry requirements, the radial problem on an annulus, and the radial problem on a ball with a Hardy potential term of critical Hardy exponent. These generalize results obtained by Sugimura \cite{Sugimura94}. 

\end{abstract}

\section{Introduction}\label{sec1}

This paper is concerned with the multiplicity of solutions to three polyharmonic Dirichlet problems. Common to all three is that the order the equations also equals the dimension of the domain. In the notation below, $d=2l$ where $d$ is the dimension, and $l$ is the power of the Laplacian. Related to this fact, they all concern odd nonlinearities of exponential growth. They also include perturbations which aren't odd. First we seek weak solutions to
\[
\tag{P}
\begin{cases}
(-\Delta)^lu = g(x,u) + \varphi(x) \quad \mbox{in $\Omega$} 
\\
\\
\bigg(\frac{\partial}{\partial \nu}\bigg)^ju\Bigg|_{\partial\Omega} = 0, \quad j = 0, \dots, l-1
\end{cases}
\] where the domain $\Omega \subset\subset \mathbb{R}^{2l}$ has a smooth boundary. Here $\varphi \in L^2(\Omega)$, and $g(x,u)$ is an odd-in-$u$, exponential nonlinearity satisfying conditions (g1)-(g5) given below. A typical example to keep in mind would be $g(x,u)=ue^{|u|^{\alpha}}$, with $0<\alpha<1$. Next in the case that $\Omega= A_{R_0}^R := \{x\in \mathbb{R}^{2l}: R_0 < |x| < R \}$ is an annulus, with $R_0>0$ and $R<+\infty$ we seek weak, radial solutions to
\[
\tag{R}
\begin{cases}
(-\Delta)^lu = 2ue^{u^2} + \varphi(x,u) \quad \mbox{in $\Omega$} 
\\
\\
\bigg(\frac{\partial}{\partial \nu}\bigg)^ju\Bigg|_{\partial\Omega} = 0, \quad j = 0, \dots, l-1
\end{cases}
\]
where $\varphi(x,u)=\varphi(|x|,u)$ is not odd in $u$. Lastly we seek weak, radial solutions to 
\[
\tag{H}
\begin{cases}
(-\Delta)^lu + \bar{b}|x|^{-2l}u= g(x,u) + \varphi(x) \quad \mbox{in $B_R(0)$} 
\\
\\
\bigg(\frac{\partial}{\partial \nu}\bigg)^ju\Bigg|_{\partial B_R(0)} = 0, \quad j = 0, \dots, l-1
\end{cases}
\] where $B_R(0)\subset \mathbb{R}^{2l}$ denotes the ball of radius $R$ centered at the origin, and $\bar{b}>0$ is a constant. Here $\varphi(x)= \varphi(|x|)$, $g(x,u)=g(|x|,u)$, and $g$ satisfies conditions (g1)-(g5) given below. We note these are not the most general possible forms such equations can take for the methods to be applicable. We have chosen them for simplicity and because they are instructive.

Many papers have been written on the existence and multiplicity of solutions for second order, nonlinear, elliptic problems, primarily by means of variational methods. The archetype has been the Dirichlet problem
\[
\begin{cases}
-\Delta u = u|u|^{p-2} + \varphi(x) \quad \mbox{in $\Omega$} \\
u|_{\partial \Omega} = 0
\end{cases}
\]
In its simplest form, the open domain $\Omega\subset \mathbb{R}^d$ is bounded, with a smooth boundary, and the perturbation $\varphi(x) \in L^2(\Omega)$. The exponent of the nonlinearity is such that if $d\geq 3$, then $2 < p < \frac{2d}{d - 2}$, while if $d=2$, then $2 < p < \infty$. These restrictions on $p$ come when one makes use of the Sobolev embeddings $W_0^{1,2}(\Omega) \hookrightarrow L^p(\Omega)$ and their compactness. If $\varphi \equiv 0$, the above equation possesses a $\mathbb{Z}_2$-symmetry with respect to the group of reflections in Sobolev space. In this case the Symmetric Mountain Pass theorem of Ambrosetti and Rabinowitz guarantees the existence of an unbounded sequence of critical values of the functional associated with the variational formulation of the problem. Such methods can be applied to a great variety of nonlinear problems invariant under compact groups of symmetries (see \cite{AmbrosettiRabinowitz1}, \cite{Rabinowitz86}, \cite{Struwe08} and references therein). However, this brings up the question of what exactly happens to this multitude of critical values when the symmetry of the problem is broken by some non-equivariant perturbation. There are no satisfactory general answers to this question yet. Methods to deal with this problem in certain cases appeared first in the papers \cite{Ba-Be1},\cite{Ba-Be2}, and \cite{Stuwe80}. The general variational principle employed in these works was later formulated by Rabinowitz in \cite{Rabinowitz82} (see also \cite{Rabinowitz86}). Roughly speaking, the idea is to estimate the spacing between consecutive symmetric mountain pass levels of the unperturbed functional, and then compare this spacing to the effect of the perturbation. Whenever the perturbation is not sufficient to eliminate this spacing, then the variational principle formulated by Rabinowitz guarantees the existence of a critical value to the perturbed functional. The first methods for estimating the spacing (or more practically the growth rate) of the symmetric mountain pass levels was based on the Weyl asymptotics for the Dirichlet eigenvalues of the Laplacian. A more refined approach, that could be more tailored for a specific problem, came in the papers of Bahri-Lions \cite{Ba-Lions}, and Tanaka \cite{Tanaka89}. Based on Morse theory, these works make use of an estimate for the number of non-positive eigenvalues of Schr\"odinger operators further described below. However, even for linear perturbations of the functional, (like the one coming from $\varphi$), the value of $p$ needs to be further restricted to $ 2 < p < \frac{2d-2}{d-2}$. To improve the range of $p$ with these methods one must weaken the perturbation. It is still a central open question of exactly how necessary is this trade-off.

Many sorts of perturbations other than a non-homogenous term $\varphi$ are of interest. A natural one is to consider the problem of an unperturbed equation, itself formally invariant, but with a non-homogenous boundary condition $u|_{\partial \Omega} = u_0 \neq 0$, which destroys the evenness of the problem. This time however, the perturbation is of much higher order, directly entering into the nonlinearity. To deal with such complications Bolle \cite{Bolle99} developed a new approach to perturbation theory of minimax levels. Similar in spirit to the earlier approach, but considerably more streamlined, the new approach considers the perturbed functional $I$ as the endpoint of a continuous path of functionals $I_{\theta}, \theta \in [0,1]$ which starts at the unperturbed functional, denoted $I_0$. Bolle's general theorem explains quantitatively how far apart two consecutive mountain pass levels of the unperturbed functional need to be for a critical level to persist to $\theta=1$. Roughly speaking, it's not the size of the perturbation at general points that determines this, but the size of $\frac{\partial}{\partial \theta}I_{\theta}(u)$ at the critical points $u$ of $I_{\theta}$. This can certainly be helpful because these $u$ satisfy the corresponding Euler-Lagrange equation. Furthermore, it becomes clearer how the size of the perturbation, as a functional in $u$, enters into the problem. This makes it easier to consider perturbations other than simple non-homogenous terms like $\varphi(x)$. This approach is further developed and applied in a number of problems in \cite{Bo-Gh-Te}.

In the two dimensional case the proper Sobolev embedding is into an Orlicz space given by an exponential $N$-function, see \cite{Adams75}. The maximal growth rate of the nonlinearity for which a variational treatment of the problem is possible is like $e^{Ku^2}$. This is related to the optimality of the Moser-Trudinger inequality \cite{Pohozaev65}, \cite{Trudinger67}, \cite{Moser71}. A typical problem now is
\[
\begin{cases}
-\Delta u = g(x,u) + \varphi(x) \quad \mbox{in $\Omega$} \\
u|_{\partial \Omega} = 0
\end{cases}
\]where $g(x,-t) = -g(x,t)$ is of exponential growth. To guarantee the convergence of general Palais-Smale sequences in $H_0^l(\Omega)$, $g$ is taken subcritical, which in this case means that $g(x,t)$ is of order of growth strictly below any positive power of $e^{t^2}$. That is
\[
\lim_{t\to \infty} \frac{|g(x,t)|}{e^{\beta t^2}} = 0 \quad \text{for all $\beta>0$.}
\]
In \cite{Sugimura94} Sugimura proved that the perturbed symmetric problem above has an infinite number of solutions if the nonlinear term $g(x,t)$ has growth like $e^{|t|^\alpha}$, where $0<\alpha<1/2$. One of the key points in this paper comes when applying the Morse index approach of \cite{Ba-Lions} and \cite{Tanaka89}. At that stage one typically applies estimates for the number of non-positive eigenvalues (aka, ``bound states") of Schr\"odinger operators. Previous results for the problem involving the power nonlinearity $u|u|^{p-2}$ had made use of a famous estimate from mathematical physics known as the CLR inequality, discovered independently by Rozenblum, Lieb, and Cwikel \cite{Rozen72},\cite{Rozen76}, \cite{Lieb76}, \cite{Cwikel77}. To get his result Sugimura proved a 2-dimensional version of this estimate. Later, using the foundation laid by Sugimura, Tarsi \cite{Tarsi08} streamlined these results by using Bolle's approach.

One advantage of Rozenblum's proof of the CLR estimate is that it automatically applies to higher order Schr\"odinger operators (which is the original form in which Rozenblum stated his result), where the Laplacian is replaced by the poly-Laplacian. Using this fact, it was proved in \cite{La-Mu-Sq} that the problem
\[
\begin{cases}
(-\Delta)^lu = u|u|^{p-2} + \varphi(x) \quad \mbox{in $\Omega$} 
\\
\\
\bigg(\frac{\partial}{\partial \nu}\bigg)^ju\Bigg|_{\partial\Omega} = \phi_j, \quad j = 0, \dots, l-1
\end{cases}
\]where $\Omega \subset\subset \mathbb{R}^{d}$, $d>2l$ has infinitely many solutions for $p$ suitably restricted.

This paper is concerned with appropriate eigenvalue estimates, and their application to problems (P), (R), and (H), which are polyharmonic, in the critical Sobolev dimension, and with exponential nonlinearities. As it is common, stronger results are proved in radially symmetric settings. For simplicity in problems (P) and (H) we consider perturbations that, in the equation, don't depend on $u$. However, more general perturbations can be considered as in \cite{Tarsi08}.

In problems (P) and (H), we make the following assumptions on the symmetric term $g$
\\
\begin{enumerate}

\item[(g1)] $g \in C(\bar{\Omega}\times \mathbb{R}, \mathbb{R})$

\item[(g2)] Given any constant $\sigma > 0$, $\exists$ constant $A_{\sigma}>0$, such that
\[
|g(x,t)| \leq A_{\sigma}e^{\sigma t^2} \quad \forall (x,t)\in \bar{\Omega}\times\mathbb{R}
\]

\item[(g3)] There are constants $\mu>0$ and $r_0\geq 0$ such that
\[
0 < G(x,t)\ln G(x,t) \leq \mu tg(x,t)
\] for $x\in \bar{\Omega}$, $|t|\geq r_0$.

\item[(g4)] $g(x,-t)=-g(x,t)$ for $(x,t)\in \bar{\Omega}\times\mathbb{R}$

\item[(g5)] There exists $0<\alpha_1 \leq \alpha_2 <1$, and $A_1, A_2, B_1$ such that
\[
A_1e^{|t|^{\alpha_1}} - B_1 \leq G(x,t) \leq A_2e^{|t|^{\alpha_2}}
\]for $(x,t)\in \bar{\Omega}\times \mathbb{R}$.
\end{enumerate}
Our first result is 
\begin{thm}\label{thm:1}
Suppose that g satisfies conditions (g1)-(g5). Then if $2/\alpha_2 - 2 > 1/\alpha_1$, problem (P) has an unbounded sequence of solutions.
\end{thm}
As a prototypical example we may take $g(x,u):= ue^{|u|^{\alpha}}$. Then the above theorem asserts that for $0<\alpha<1/2$ problem (P) has an unbounded sequence of solutions.  If we impose radial assumptions on the problem we can improve the result. In problem (R), let $\Phi(x,u):= \int_{0}^{u}\varphi(x,t)dt$
\begin{thm}\label{thm:02}
Suppose that there exists a $\beta \in (0,1)$ such that $|\Phi(x,t)|+|\varphi(x,t)t|\leq C(t^2e^{t^2})^{\beta}$ for sufficiently large values of $|t|$. If $2l>\frac{1}{1-\beta}$ then problem (R) has an unbounded sequence of radial solutions.
\end{thm}
\begin{thm}\label{thm:2}
Suppose that $\Omega=B_R(0)$, the open ball centered at the origin with a finite radius $R$. Suppose that, in addition to conditions (g1)-(g5), we have $g(x,u)=g(|x|,u)$ and $\varphi(x)=\varphi(|x|)$. Then if $2/\alpha_2 - 1 > 1/\alpha_1$, problem (H) has an unbounded sequence of radial solutions.
\end{thm}
These theorems are proven using upper estimates on the number of non-positive eigenvalues (a.k.a. bound states) of polyharmonic Schr\"odinger operators. When no radial symmetry is imposed, the corresponding spectral inequality takes the following form:
\begin{proposition}\label{prop:1}
Let $\mathcal{B}(t):= (|t|+1)\ln(|t|+1)-|t|$ be an $N$-function, and $L_{\mathcal{B}}(\Omega)$ be the corresponding Orlicz space. Let $V(x)\in L_{\mathcal{B}}(\Omega)$, and on $L^2(\Omega)$ consider the unbounded linear operator $(-\Delta)^l -V(x)$. Denote by $N_-((-\Delta)^l-V(x))$ the number of its non-positive eigenvalues. Then there exists a constant $C = C(l,\Omega)$ such that
\[
N_-((-\Delta)^l-V(x)) \leq C||V||_{L_{\mathcal{B}(\Omega)}}
\]
\end{proposition}
Simplified for our purposes, this result is due to Solomyak \cite{Solomyak94}. We will not use Proposition 1 but rather a corollary of it. As is well-known in the theory of Orlicz spaces $||V||_{L_{\mathcal{B}}(\Omega)} \leq \int_{\Omega}\mathcal{B}(V(x))dx+1$. See (9.12) in \cite{Kras-Rutic}.
\begin{coro}\label{cor:1}
$N_-((-\Delta)^l-V(x)) \leq C\int_{\Omega}\mathcal{B}(V(x))dx+C$.
\end{coro}
This generalizes the eigenvalue estimate originally obtained by Sugimura \cite{Sugimura94}. It is unknown to the author if the finer estimate in Proposition \ref{prop:1} can yield stronger results in this case. To answer such a question a much finer analysis of the critical sequence given by Tanaka's theorem below (or a possible alternative of it) seems necessary. In the radial setting on an annulus the main eigenvalue estimate we prove is the following 
\begin{proposition}\label{prop:2}
Let $\Omega$ be an annulus of outer radius $R<+\infty$ and inner radius $R_0>0$. Let $V(x)=V(|x|)$. On the space $L_r^2(\Omega)$, of radially symmetric, square integrable functions, consider the unbounded linear operator $(-\Delta)^l -V(x)$. Denote by $N_-((-\Delta)^l-V(x))$ the number of its non-positive eigenvalues. There exists a constant $C = C(l,R_0,R)$ such that
\[
\big[N_-((-\Delta)^l-V(x))\big]^{2l} \leq C\int_{\Omega}V^+(x)\bigg[1+ \log\bigg(\frac{R}{|x|}\bigg)\bigg]^{2i}dx
\]where $i = 1/2$ when $l=1$, and $i=1$ when $l>2$.
\end{proposition}

\section{Preliminaries: General outline of the perturbation method}\label{sec2}

{\it 2.1 The symmetric mountain pass sequence and the variational principle of Bolle}\\

Here we recall the perturbation method of Bolle. See \cite{Bolle99} and \cite{Bo-Gh-Te}. Let $E$ be a Hilbert Space, and let ${e_k}$ be a basis for $E$. Decompose $E$ as 
\begin{equation}\label{eq:1}
E = \overline{\cup^{\infty}_{k=0}E_k}, \quad \text{where} \quad E_k = E_{k-1} \oplus \mathbb{R}e_k
\end{equation} with $E_0=\{0\}$. For a given increasing sequence of real numbers $R_k > 0$ set
\[
\Gamma_k := \{g \in C(E_k,E):\text{$g$ is odd and $g(u)=u$ if $u\in E_k$  $||u|| \geq R_k$}\}
\] The precise values of $R_k$ that we will use are defined later. For an even functional $I_0(u)$ on $E$ set
\[
c_k := \inf_{g\in \Gamma_k}\sup_{u \in g(E_k)}I_0(u)
\]
Under typical assumptions on $I_0(u)$, the Symmetric Mountain Pass Theorem of Ambrosetti and Rabinowitz shows that, for appropriate values of $R_k$ the minimax values $c_k$ are an unbounded sequence of critical values of $I_0(u)$.

Consider two continuous functions $f_1,f_2:[0,1]\times \mathbb{R} \to \mathbb{R}$ which are Lipschitz continuous with respect to the second variable, and $f_1\leq f_2$. Define the associated flows
\[
\begin{cases}
\psi_i(0,s) = s\\
\frac{\partial}{\partial\theta}\psi_i(\theta,s) = f_i(\theta, \psi_i(\theta,s))
\end{cases}
\]$\psi_1$ and $\psi_2$ are continuous in both variables, and non-decreasing in $s$. Moreover, by the comparison theorem for ODEs, since $f_1\leq f_2$, we have $\psi_1 \leq \psi_2$. Also, denote
\[
\bar{f_i}(s) := \sup_{\theta \in [0,1]}|f_i(\theta,s)|
\]We will apply the main theorem of \cite{Bo-Gh-Te}:

\begin{thm}[Bolle, Ghoussoub, Tehrani \cite{Bo-Gh-Te}]\label{thm:3}
Let $E$ be a Hilbert space and $I:[0,1]\times E \to \mathbb{R}$ be a $C^1$ functional satisfying the following conditions
\\
\\
(H1) $I$ satisfies the following analogue of the Palais-Smale Condition: For a sequence $\{(\theta_n,u_n)\}_{n\in \mathbb{N}} $ in $[0,1]\times E$ such that $||I_{\theta_n}'(u_n)||_{E*} \to 0$ and $|I_{\theta_n}(u_n)| \leq C$ there is a subsequence of it strongly converging in $[0,1]\times E$.
\\
\\
(H2) For all $b>0 \:\: \exists$ a positive constant $C(b)$ such that
\[
|I_{\theta}(u)| \leq b \quad \text{implies} \quad \bigg|\frac{\partial}{\partial\theta}I_{\theta}(u)\bigg| \leq C(b)(||I_{\theta}'(u)||+1)(||u||+1)
\]
(H3) There exist two continuous functions $f_1,f_2:[0,1]\times\mathbb{R}\to\mathbb{R}$, with $f_1 \leq f_2$, that are Lipschitz continuous relative to the second variable, and such that, for all critical points u of $I_{\theta}$
\[
f_1(\theta,I_{\theta}(u)) \leq \frac{\partial}{\partial \theta}I_{\theta}(u) \leq f_2(\theta, I_{\theta}(u))
\]
(H4) $I_0$ is even and for any finite dimensional subspace $W$ of $E$ we have
\[
\sup_{\theta \in [0,1]}I_{\theta}(y) \to -\infty \:\: \text{as $||y|| \to \infty$ for $y\in W$}.
\]Then there is a $K>0$ such that for every $n$ only one of the two possibilities below holds:
\\
(1) Either $I_1$ has a critical level $\bar{c}_n$ with
\[
\psi_2(1,c_n) < \psi_1(1,c_{n+1}) \leq \bar{c}_n,
\] 
(2) Or $c_{n+1}-c_n \leq K(\bar{f_1}(c_{n+1})+\bar{f_2}(c_n)+1)$.

\end{thm}

The values $R_k$ that we will use are defined as follows: by hypothesis (H4) in the theorem above we can find an $R_k>0$ such that $\sup_{\theta \in [0,1]}I_{\theta}(u)<0$ for all $u\in E_k$ with $||u||\geq R_k$.

The theorem above shows that, for each $n$, two outcomes are possible. If the second outcome holds for all sufficiently large $n$ then one can derive an upper bound on the sequence $c_n$. In applications, one shows that the first possibility holds for infinitely many $n$ by obtaining a lower bound for $c_n$ that contradicts this alleged upper bound. Because $c_n\to+\infty$ and because in typical applications $\psi_2(1,s)\to+\infty$ as $s \to +\infty$ we get that $\bar{c}_n \to +\infty$. We also note that the above theorem was originally proven for $I$ in $C^2(E)$, by using the gradient flow of $I$. However, as is typical, it suffices for the proof to use only a pseudo-gradient flow in the sense of Palais. This requires that $I$ is only in $C^1(E)$. See \cite{Clapp-Ding-Hernandez}.\\
\\

\section{Problem (P), the case of an unrestricted domain}\label{sec3}

{\it 3.1 The variational setup}\\

Let $\Omega \subset \subset \mathbb{R}^{d}$ be a smooth domain, and where ${d=2l}$. In the Hilbert space $L^2(\Omega)$ we consider the dense subspace
\[
H^l_0(\Omega) := \text{the completion of $C^{\infty}_0(\Omega)$ in the following norm}
\]

\begin{equation}\label{eq:2}
||u||_{H^l_0(\Omega)} = \begin{cases}
||\Delta^k u||_{L^2(\Omega)} & \mbox{if $l$=$2k$}\\
||\nabla(\Delta^k u)||_{L^2(\Omega)} & \mbox{if $l$ = $2k+1$}
\end{cases}
\end{equation} 

Typically on the space $H^l_0(\Omega)$ we use the norm $||D^lu||_{L^2(\Omega)}$, after taking into account Poincare's Inequality for the lower order terms. However, this norm is equivalent to $||u||_{H^l_0(\Omega)}$ on $C^{\infty}_0(\Omega)$ by integration by parts. When convenient, we shall denote the norm $||u||_{H^l_0(\Omega)}$ by $||u||$. As shorthand, on $H^l_0(\Omega)$, define the $l$th power of the gradient as 

\begin{equation}\label{eq:3}
\nabla^lu = \begin{cases}
\Delta^ku & \mbox{if $l=2k$}\\
\nabla (\Delta^ku) & \mbox{if $l=2k+1$}
\end{cases}
\end{equation}

The variational setup for our problem is as follows. On the space $H^l_0(\Omega)$ we consider the functional
\begin{equation}\label{eq:4}
I_1 := \frac{1}{2}\int_{\Omega}|\nabla^lu|^2dx - \int_{\Omega} G(x,u)dx - \int_{\Omega}\varphi u dx
\end{equation} This is the functional whose critical points correspond to generalized solutions of the boundary value problem
\[
\label{P}
\tag{P}
\begin{cases}
(-\Delta)^lu = g(x,u) + \varphi(x), \quad \mbox{$x\in \Omega$} 
\\
\\
\Bigg(\frac{\partial}{\partial \nu}\Bigg)^ju\Bigg|_{\partial\Omega} = 0, \quad j = 0, \dots, l-1
\end{cases}
\]where $G(x,u) = \int_{0}^{u}g(x,s)ds$, in the space $H^l_0(\Omega)$. See \cite{Gazzola-Grunau-Sweers}.

For $I_1(u)$ to be well-defined on all of $H^l_0(\Omega)$, and for such a variational treatment to be viable, we must restrict the growth rate with respect to $u$ of the nonlinearity $G(x,u)$. The maximal growth for which a variational treatment is allowed comes from the optimality of Adam's generalization of the Moser-Trudinger inequality. Namely, on the space $W^{l,\frac{d}{l}}_0(\Omega)$, $1\leq l<d$, Adams showed

\begin{equation}\label{eq:5}
\sup\limits_{\substack{u \in W^{l,\frac{d}{l}}_0(\Omega) \\ ||\nabla^lu||_\frac{d}{l} \leq 1}}  \int_{\Omega} e^{\beta|u|^{\frac{d}{d-l}}}dx 
\begin{cases}
\leq C|\Omega|, \quad &\text{if $\beta \leq \beta(d,l)$} \\
= +\infty, \quad &\text{if $\beta > \beta(d,l)$}
\end{cases}
\end{equation} where $\beta(d,l)$ is given explicitly. See \cite{DRAdams88}. In our case, $d=2l$ so the exponent $\frac{d}{d-l}$ equals 2. In this case, $\beta(2l,l) = l!(4\pi)^l$. For a variational treatment to be possible in $H_0^l(\Omega)$, $G(x,u)$ can't grow faster than $e^{Ku^2}$ for all $K$.\\

Conditions (g1)-(g5) imply that $I_1(u)$ is a $C^1$ functional on $H^l_0(\Omega)$. For the corresponding path of functionals we simply consider

\begin{equation}\label{eq:6}
I_{\theta}(u) := \frac{1}{2}\int_{\Omega}|\nabla^lu|^2dx - \int_{\Omega} G(x,u)dx - \theta\int_{\Omega}\varphi u dx 
\end{equation} where $\theta \in [0,1]$, and for which $I_0$ is even.
\\
\\
\\
{\it 3.2 Bolle's Requirements}\\
\\
{\bf(H1)} $I$ satisfies the following analogue of the Palais-Smale Condition: For a sequence $\{(\theta_n,u_n)\}_{n\in \mathbb{N}} $ in $[0,1]\times E$ such that $||I_{\theta_n}'(u_n)||_{E*} \to 0$ and $|I_{\theta_n}(u_n)| \leq C$ there is a subsequence converging strongly in $[0,1]\times E$.
\\
\\
{\bf Proof}: Let $\{(\theta_n,u_n)\}_{n\in \mathbb{N}}$ be such a sequence. Then, after taking a subsequences, we can find constants $C_0$ and $\theta_0$ such that
\begin{equation}\label{eq:7}
I_{\theta_n}(u_n) =  \frac{1}{2}\int_{\Omega}|\nabla^lu_n|^2dx - \int_{\Omega} G(x,u_n)dx - \theta_n\int_{\Omega}\varphi u_n dx \to C_0 
\end{equation} where $\theta_n\to \theta_0$, and
\begin{equation}\label{eq:8}
\bigg| \int \nabla^lu_n \cdot \nabla^lv -g(x,u_n)v - \theta_n\varphi(x)v \bigg| \leq \epsilon_n||v||_{H^l_0(\Omega)}
\end{equation}for all $v\in H^l_0(\Omega)$, where $\epsilon_n \to 0$ as $n \to \infty$. Choosing $v= u_n$ in (8) we get, for $\bar{\mu} > 2$
\begin{equation}\tag{7'}\label{eq:7'}
\frac{\bar{\mu}}{2}\int_{\Omega}|\nabla^lu_n|^2dx - \bar{\mu}\int_{\Omega} G(x,u_n)dx - \bar{\mu}\theta_n\int_{\Omega}\varphi(x) u_n dx \leq C_1
\end{equation} and
\begin{equation}\label{eq:9}
 -||u_n||^2_{H^l_0(\Omega)} + \int_\Omega g(x,u_n)u_ndx + \theta_n\int_\Omega\varphi(x)u_ndx   \leq \epsilon_n||u_n||_{H^l_0(\Omega)}
\end{equation} Adding (\ref{eq:7'}) and (\ref{eq:9}) gives
\begin{eqnarray}\label{eq:10}
\bigg(\frac{\bar{\mu}}{2}-1 \bigg)||u_n||^2 + \int_{\Omega} \bigg(g(x,u_n)u_n - \bar{\mu}G(x,u_n)\bigg)dx + (1-\bar{\mu})\theta_n\int_{\Omega}\varphi(x)u_ndx\nonumber\\ \leq C' + \epsilon_n||u_n||
\end{eqnarray}By assumption (g3)
\[
\int_{\Omega} \bigg(g(x,u_n)u_n - \bar{\mu}G(x,u_n)\bigg)dx \geq -C_2
\]Hence (\ref{eq:10}) gives
\begin{eqnarray}
\bigg(\frac{\bar{\mu}}{2}-1 \bigg)||u_n||^2 &\leq& C_3 + \epsilon_n||u_n|| + C_4\int|\varphi(x)u_n|dx\nonumber\\
&\leq& C_3 + \epsilon_n||u_n|| + C_5||u_n||_{L^2(\Omega)}\nonumber
\end{eqnarray}
By the generalized Poincar\'e Inequality if $u\in H^l_0(\Omega)$ then there exists a constant $C_6>0$ such that
\[
||u||_{L^2(\Omega)}\leq C_6||u||_{H^l_0(\Omega)}
\] So we get
\[
\bigg(\frac{\bar{\mu}}{2}-1 \bigg)||u_n||^2 \leq C_3 + C_7||u_n||
\]Thus
\begin{equation}\label{eq:11}
||u_n|| \leq K
\end{equation}
Having proven the boundedness of Palais-Smale sequences, we will show that they are pre-compact by proving that $I_{\theta}'(u)(\cdot)$ has the form $L(u)(\cdot) + K(u)(\cdot)$ where $L:H^l_0(\Omega) \to H^{-l}(\Omega)$ is an isomorphism, and $K:H^l_0(\Omega) \to H^{-l}(\Omega)$ is compact.  Although this isn't entirely necessary, and a shorter proof which doesn't rely explicitly on this fact is possible. However, the fact that $I_{\theta}'(u)$ has this form will be needed later to apply Tanaka's Theorem, hence we prove it now. Now
\begin{equation}\label{eq:12}
I_{\theta_0}'(u)(\cdot) = \langle u,\cdot\rangle_{H^l_0(\Omega)} - \langle g(x,u),\cdot\rangle_{L^2(\Omega)}-\theta_0 \langle \varphi(x),\cdot\rangle_{L^2(\Omega)}
\end{equation}Clearly $L:H^l_0(\Omega) \to H^{-l}(\Omega): u \mapsto \langle u,\cdot \rangle_{H^l_0(\Omega)} $ is the Riesz map, hence a Hilbert space isomorphism. Clearly the map $K_1: H^l_0(\Omega) \to H^{-l}(\Omega): u \mapsto \langle \varphi(x),\cdot\rangle_{L^2(\Omega)}$ is compact because it's a constant map. To show that $K_2:H^l_0(\Omega) \to H^{-l}(\Omega): u \mapsto \langle g(x,u),\cdot\rangle_{L^2(\Omega)}$ is compact it suffices to show that if $\{ u_n \} \subset H^l_0(\Omega)$ is bounded then, up to a subsequence, $g(x,u_n)$ converges in $L^2(\Omega)$. WLOG we may assume, after taking a subsequence, that
\\
\\
$||u_n|| \leq K$
\\
$u_n \rightharpoonup u$ weakly in $H^l_0(\Omega)$
\\
$u_n \to u$ strongly in $L^p(\Omega)$, $p \geq 1$
\\
$u_n(x) \to u(x)$ a.e in $\Omega$
\\
\\
Now since $g$ has subcritical growth in $u$ by (g2), we can find $C_K>0$ such that
\begin{equation}\label{eq:13}
|g(x,t)| \leq C_K \exp \bigg( \frac{\beta(2l,l)}{K^2}t^2 \bigg)
\end{equation} where $\beta(2l,l)=l!(4\pi)^l$ is the optimal constant in Adam's inequality. We apply Adam's inequality:
\begin{eqnarray}
||g(x,u_n)||_{L^2}^2 &\leq& C_K\int_{\Omega}\exp\bigg( \frac{\beta(2l,l)}{K^2}|u_n|^2 \bigg)dx\nonumber\\
&\leq& C_K\int_{\Omega}\exp\bigg( \frac{\beta(2l,l)}{||u_n||^2}|u_n|^2 \bigg)dx\nonumber\\
&\leq C_K'\nonumber&
\end{eqnarray} Similarly we have
\[
\int_{\Omega}|g(x,u_n)|^2|u_n|dx \leq C_K''
\]To obtain the required result we use the following lemma

\begin{lemma}
Let $\{u_n\}$ be a convergent sequence of functions in $L^2(\Omega)$, with $u_n(x) \to u(x)$ a.e. Assume that $g(x,u_n)$ and $g(x,u)$ are also in $L^2(\Omega)$ with $g(x,t)$ continuous in $t$ uniformly in $x$. If  
\[
\int_{\Omega}|g(x,u_n)|^2|u_n|dx \leq C_8
\]then $g(x,u_n)$ converges in $L^2(\Omega)$ to $g(x,u)$.
\end{lemma}
{\bf Proof}: Note that since $g(x,t)$ is continuous in $t$ and $u_n(x)\to u(x)$ a.e. then $g(x,u_n(x))\to g(x,u(x))$ a.e. We have that
\begin{eqnarray}
|g(x,u_n(x))-g(x,u(x))|^2 &\leq& \big[ |g(x,u_n(x))|+ |g(x,u(x))| \big]^2\nonumber\\
&\leq& 2|g(x,u_n(x))|^2+ 2|g(x,u(x))|^2\nonumber
\end{eqnarray}If we assume $||g(x,u_n)||_{L^2}\to||g(x,u)||_{L^2}$ then 
\[
\int_{\Omega}2|g(x,u_n(x))|^2+ 2|g(x,u(x))|^2dx \to 4\int_{\Omega}|g(x,u(x))|^2dx
\]Also $|g(x,u_n(x))-g(x,u(x))|\to 0$ a.e. So we can apply the generalized Lebesgue Dominated Convergence Theorem and get that $\int|g(x,u_n(x))-g(x,u(x))|^2dx\to 0$, which is the required result. So it suffices to prove that $\int |g(x,u_n)|^2\to \int|g(x,u)|^2dx$. Let $f(x,t):= g(x,t)^2$

Since $f(x,u(x)) \in L^1(\Omega)$ it follows that for a given $\epsilon>0$ there is a $\delta>0$ such that
\begin{equation}\label{eq:14}
\int_A f(x,u(x))dx \leq \epsilon \quad \text{if $|A|\leq \delta$}
\end{equation} for all measurable subsets $A\subseteq \Omega$. Next using the fact that $u\in L^1(\Omega)$ we find $M_1>0$ such that
\begin{equation}\label{eq:15}
|\{ x\in\Omega: |u(x)|\geq M_1 \}|\leq \delta
\end{equation}
Let $M:= \max\{M_1,C_8/\epsilon\}$. We write
\begin{equation}\label{eq:16}
\bigg|\int f(x,u_n(x))dx  -\int f(x,u(x))dx\bigg| \leq I_1 + I_2 + I_3
\end{equation} and estimate each integral separately:
\begin{eqnarray}
I_1 \equiv \int_{|u_n(x)|\geq M}f(x,u_n(x))dx &=& \int_{|u_n(x)|\geq M}\frac{|g(x,u_n(x))|^2}{|u_n(x)|}|u_n(x)|dx\nonumber\\
&\leq& \frac{C_8}{M} \leq \epsilon\nonumber
\end{eqnarray}By the choices we have made above (\ref{eq:14}) and (\ref{eq:15}) imply that
\[
I_3 \equiv \int_{|u(x)|\geq M}f(x,u(x))dx \leq \epsilon
\]Next we claim that 
\[
I_2 \equiv \bigg| \int_{|u_n(x)|< M}f(x,u_n(x))dx - \int_{|u(x)|< M}f(x,u(x))dx \bigg| \to 0
\] as $n\to \infty$.
Indeed, $h_n(x):= f(x,u_n(x))\chi_{|u_n|<M} - f(x,u(x))\chi_{|u|<M}$ tends to 0 a.e. in $\Omega$. Moreover $|h_n(x)|\leq |f(x,u(x))|$ if $|u_n(x)|\geq M$ and $|h_n(x)| \leq C + f(x,u(x))$ if $|u_n(x)|<M$. So $I_2 \to 0$ as $n\to \infty$ by the Lebesgue Dominated Convergence Theorem. $\blacksquare$
\\
\\
Thus $I'$ has the stated form and (H1) is satisfied.
\\
\\
\\
{\bf (H2)} $\: \: \:$ Here $\frac{\partial}{\partial\theta}I_{\theta}(u) = -\int_{\Omega}\varphi(x) u(x)dx$ is bounded in absolute value by $||\varphi||_{L^2(\Omega)}||u||_{L^2(\Omega)}$. By the Generalized Poincar\'e inequality this is bounded by $C_{\varphi}||u||$.
\\
\\
\\
{\bf (H3) Determining $f_1,f_2$}
\\
\begin{lemma}
There exists a constant $C_9>0$ such that if $u\in H^l_0(\Omega)$ is a critical point of $I_{\theta}$ then
\begin{equation}\label{eq:17}
\bigg| \frac{\partial}{\partial\theta}I_{\theta}(u)\bigg| \leq C_9 [\ln(|I_{\theta}(u)|+1) ]^{1/\alpha_1} +C_9
\end{equation}
\end{lemma}
{\bf Proof}: From (H2) above we have that
\[
\bigg|\frac{\partial}{\partial\theta}I_\theta(u)\bigg|\leq C_{10}||u||_{L^2(\Omega)}.
\]
So it suffices to estimate $||u||_{L^2(\Omega)}$. Assume $I_{\theta}'(u)=0$. Then
\begin{eqnarray}\label{eq:18}
I_{\theta}(u) &=&I_{\theta}(u) - \frac{1}{2}\langle I_{\theta}'(u),u \rangle \nonumber\\
&=&\frac{1}{2}\int_{\Omega}|\nabla^lu|^2dx - \int_{\Omega}G(x,u)dx -\theta\int_{\Omega}\varphi udx \nonumber\\
& &- \frac{1}{2}\int_{\Omega}|\nabla^lu|^2dx + \int_{\Omega}\frac{1}{2} g(x,u)udx +\frac{\theta}{2}\int_{\theta}\varphi udx\nonumber\\\nonumber\\
&=& \int_{\Omega}\bigg(\frac{1}{2} g(x,u)u-G(x,u)\bigg)dx - \frac{\theta}{2}\int_{\theta}\varphi udx
\end{eqnarray}
We now apply condition (g3). When $|u(x)|>r_0$ we bound
\begin{eqnarray}
\frac{1}{2}g(x,u)u -G(x,u) &\geq& \frac{1}{2\mu}G(x,u)\ln[G(x,u)]-G(x,u)\nonumber\\
&\geq& \bigg( \frac{1}{2\mu} -\epsilon \bigg)G(x,u)\ln[G(x,u)]
-C_{11}\nonumber
\end{eqnarray}
When $|u(x)|\leq r_0$ the expression $\frac{1}{2}g(x,u)u -G(x,u)$ is bounded by a constant since $g$ and $G$ are continuous. Since $\Omega$ is of finite measure we get from the above and equation (\ref{eq:18})
\begin{equation}
I_{\theta}(u) \geq \bigg( \frac{1}{2\mu} -\epsilon \bigg)\int_{|u(x)|\geq r_0} G(x,u)\ln[G(x,u)]dx -C_{12}||u||_{L^2(\Omega)} - C_{13} \nonumber
\end{equation}
Now applying the growth condition (g5) and again the fact that $\Omega$ is of finite measure
\begin{equation}\label{eq:19}
I_{\theta}(u) \geq C_{14}\int_{\Omega}|u|^{\alpha_1}e^{|u|^{\alpha_1}}dx -C_{12}||u||_{L^2(\Omega)} -C_{13}
\end{equation}
Now observe that for $\alpha, \beta > 0$ there exists a constant $t_0=t_0(\alpha,\beta)$ such that the function $t^{\beta}e^{t^{\alpha}}$ is convex for $t\geq t_0$. We take $\alpha = \beta = \alpha_1/2$, and apply Jensen's inequality
\begin{equation}\label{eq:20}
C_{15}\int_{|u|\geq t_0^{1/2}}|u|^{\alpha_1}e^{|u|^{\alpha_1}}dx \geq \bigg\{ \frac{1}{|\Omega|}\int_{|u|\geq t_0^{1/2}}|u|^2dx \bigg\}^{\alpha_1/2} \exp\bigg[ \bigg\{ \frac{1}{|\Omega|}\int_{|u|\geq t_0^{1/2}}|u|^2dx\bigg\}^{\alpha_1/2}\bigg]
\end{equation}
Also note that 
\begin{equation}\label{eq:21}
||u||_{L^2(\Omega)}^{\alpha_1} \leq C_{16} + \bigg\{ \int_{|u|\geq t_0^{1/2}} |u|^2 dx \bigg\}^{\alpha_1/2}
\end{equation}
Hence we proceed as follows
\begin{multline}
||u||_{L^2}^{\alpha_1}\exp\bigg[ \bigg(\frac{1}{|\Omega|}\bigg)^{\alpha_1/2}||u||_{L^2}^{\alpha_1} \bigg] \quad \leq \nonumber\\ \bigg\{C_{16} + \bigg( \int_{|u|\geq t_0^{1/2}} |u|^2 dx \bigg)^{\alpha_1/2}\bigg\}\exp\bigg[ C_{16} + \bigg\{ \frac{1}{|\Omega|}\int_{|u|\geq t_0^{1/2}} |u|^2 dx \bigg\}^{\alpha_1/2} \bigg]\\\\
\leq C_{17} +  C_{18}\bigg\{ \int_{|u|\geq t_0^{1/2}} |u|^2 dx\bigg\}^{\alpha_1/2}\exp\bigg[ \bigg\{ \frac{1}{|\Omega|}\int_{|u|\geq t_0^{1/2}} |u|^2 dx \bigg\}^{\alpha_1/2} \bigg]
\end{multline}Here apply (\ref{eq:20})
\begin{eqnarray}\label{eq:22}
&\leq& C_{19} + C_{20}\int_{|u|\geq t_0^{1/2}}|u|^{\alpha_1}e^{|u|^{\alpha_1}}dx\nonumber\\
&\leq& C_{19} + C_{20}\int_{\Omega}|u|^{\alpha_1}e^{|u|^{\alpha_1}}dx
\end{eqnarray}
Inequality (\ref{eq:22}) now implies 
\begin{equation}\label{eq:23}
||u||_{L^2}\leq C_{21}\bigg\{ \ln \bigg( \int_{\Omega}|u|^{\alpha_1}e^{|u|^{\alpha_1}}dx +1 \bigg) \bigg\}^{1/\alpha_1}+C_{22}
\end{equation}
So (\ref{eq:23}) and (\ref{eq:19}) give
\begin{eqnarray}\label{eq:24}
I_{\theta}(u)&\geq& C_{23}\int_{\Omega}|u|^{\alpha_1}e^{|u|^{\alpha_1}}dx - \bigg\{ \ln\bigg( \int_{\Omega}|u|^{\alpha_1}e^{|u|^{\alpha_1}}dx+1\bigg)\bigg\}^{1/\alpha_1} -C_{24}\nonumber\\
&\geq& C_{25}\int_{\Omega}|u|^{\alpha_1}e^{|u|^{\alpha_1}}dx-C_{26}
\end{eqnarray}
So by (\ref{eq:22}) and (\ref{eq:24})
\begin{equation}\label{eq:25}
I_{\theta}(u) \geq C_{27}||u||_{L^2}^{\alpha_1}\exp\bigg[ \bigg(\frac{1}{|\Omega|}\bigg)^{\alpha_1/2}||u||_{L^2}^{\alpha_1} \bigg]-C_{28}\nonumber
\end{equation}
that is,
\begin{equation}
||u||_2 \leq C_{29} [\ln(|I_{\theta}(u)|+1) ]^{1/\alpha_1} +C_{30}
\end{equation}which proves the lemma.
$\blacksquare$
\\
\\
Thus we take 
\begin{equation}\label{eq:26}
f_i(\theta,t) = f_i(t) = (-1)^iC_9\big\{ \big[ \ln(|t|+1) \big]^{1/\alpha_1} +1\big\}
\end{equation}
\\
\\
{\bf (H4)} This condition is easily satisfied by assumption (g5), which shows that $G(x,u)$ is super-quadratic (uniformly in $x$) and tends to $+\infty$ as $|u|\to \infty$.\\
\\
\\
{\it 3.3 The Alleged Upper Bound}\\

As noted earlier, we will operate under the assumption that alternative 2) of Theorem 2 holds for sufficiently large $n$.  That is, for $n>n_0$, it'll assumed that
\begin{equation}\label{eq:27}
c_{n+1}-c_n \leq K\big[ \ln(c_{n+1})^{1/\alpha_1} + \ln(c_n)^{1/\alpha_1} + 1 \big]
\end{equation}
We will show this implies that $c_n \leq An[\ln(n)]^{1/\alpha_1}$ for sufficiently large $n$, and for some constant $A$ to be chosen appropriately. \\

Let $\gamma := 1/\alpha_1$ and let $b_n := An[\ln(n)]^{\gamma}$. First we can choose $A>0$ so large that $c_{n_0} < b_{n_0}$ where $n_0$ is large and fixed. For $n> n_0$
\begin{eqnarray}
b_{n+1}-b_{n} &=& A(n+1)[\ln(n+1)]^{\gamma} - An[\ln(n)]^{\gamma}\nonumber\\
&=& A\big[ \ln(n+\theta)^{\gamma} + \gamma\ln(n+\theta)^{\gamma-1} \big]\nonumber
\end{eqnarray}for some $\theta\in[0,1]$ by the Mean Value Theorem. Hence
\begin{equation}\label{eq:28}
b_{n+1}-b_{n} \geq \frac{A}{2}\big[ \ln(n) \big]^{\gamma}
\end{equation}
Now from the definition of $b_n$ we compute
\begin{equation}\label{eq:29}
K\big[\ln(b_{n+1})^{\gamma}+\ln(b_n)^{\gamma}+1\big] \leq C_{\gamma}K\ln(n)^{\gamma} + C_{\gamma}K\ln(A)
\end{equation} for $n>n_0$ sufficiently large. So we take $A>>2C_{\gamma}K$. Then (\ref{eq:28}) and (\ref{eq:29}) combine to give
\[
b_{n+1}-b_n > K\big[\ln(b_{n+1})^{\gamma}+\ln(b_n)^{\gamma}+1\big]
\]This is the reverse of the inequality satisfied by $c_n$. We already have $b_{n_0} \geq c_{n_0}$. Assume that $b_i>c_i$ for $i = n_0,\dots,n$. We will show that $b_{n+1}\geq c_{n+1}$:
\begin{eqnarray}\label{eqApp429}
b_{n+1}-c_{n+1} &=& b_{n+1}-b_n - (c_{n+1}-c_n) + (b_n-c_n)\nonumber\\
&\geq& b_{n+1}-b_n - (c_{n+1}-c_n)\nonumber\\
&\geq& K\big[\ln(b_{n+1})^\gamma +\ln(b_n)^\gamma+1\big]-K\big[\ln(c_{n+1})^\gamma +\ln(c_n)^\gamma+1\big]\nonumber\\
&=&K\ln(b_{n+1})^\gamma-K\ln(c_{n+1})^\gamma+ \big[K\ln(b_{n})^\gamma-K\ln(c_{n})^\gamma\big]\nonumber\\
&\geq& K\ln(b_{n+1})^\gamma-K\ln(c_{n+1})^\gamma.
\end{eqnarray}
Assume that $b_{n+1} < c_{n+1}$. Then
\begin{eqnarray}\label{eqApp430}
K\ln(b_{n+1})^\gamma - K\ln(c_{n+1})^\gamma &=& -K\gamma\int_{b_{n+1}}^{c_{n+1}}\frac{\ln(t)^{\gamma-1}}{t}dt\nonumber\\
&>& b_{n+1}-c_{n+1},
\end{eqnarray}
where we have used the fact that $-K\gamma\frac{\ln(t)^{\gamma-1}}{t}>-1$ for $t>b_{n+1}\geq b_{n_0}$, when $n_0$ is taken to be sufficiently large. This contradicts (\ref{eqApp429}), and so $c_{n+1}\leq b_{n+1}$. Thus by induction $c_n \leq b_n$ for all $n>n_0$, i.e.
\begin{equation}\label{eqApp29}
c_n \leq An\big[\ln(n)\big]^{1/\alpha_1} \quad \text{for $n>n_0$},
\end{equation}when assuming alternative 2) of Theorem \ref{thm:3}.

\section{Tanaka's Theorem and its Requirements}\label{sec4}

The reference is Theorem B in \cite{Tanaka89} (see also \cite{Ba-Lions}). The idea is that associated to each minimax value of the symmetric functional there is a sequence of critical points which are at a lower energy level, but which have a large augmented Morse index. For the moment, we are concerned with $I_0(u)$
\[
I_0(u) = \frac{1}{2}\int_{\Omega}|\nabla^lu|^2dx - \int_{\Omega} G(x,u)dx
\]

We already know from the proof of the Palais-Smale condition that $I_0'$ has the form of a compact perturbation of a Hilbert space isomorphism. Actually we will apply Tanaka's Theorem to a slightly smoother functional:
\begin{equation}\label{eq:31}
J(u):=\frac{1}{2}\int_{\Omega}|\nabla^lu|^2dx - \int_{\Omega} H(u)dx
\end{equation}where
\[
H(t) = a\exp[(t^2+1)^b]
\] where $b=\alpha_2/2$. By assumption (g5) we can choose $a>0$ so that
\[
G(x,t) \leq H(t) \quad \text{for $(x,t) \in \bar\Omega\times\mathbb{R}$}
\]Thus $I_0(u)\geq J(u)$.

$J(u)$ has a nonlinearity of subcritical and super-quadratic growth, and so all compactness properties of $I_0(u)$ also hold for $J(u)$. In particular $J'$ has the form $L+K$ where $L:H^l_0(\Omega) \to H^{-l}(\Omega)$ is an isomorphism, and $K:H^l_0(\Omega) \to H^{-l}(\Omega)$ is compact. In addition we have the following compactness conditions needed in the application of Tanaka's theorem: Let $\{ E_j \}$ be the decomposition in equation (\ref{eq:1}).
\\
\\
$(PS)_m \quad$If for some $M>0$, $\{u_j\}$ satisfies
\[
u_j\in E_m, \quad J(u_j)\leq M \quad \forall j, \quad ||(J|_{E_m})'(u_j)||_{E'_m} \to 0 \quad \text{as $j\to \infty$}
\]then $\{u_j\}$ is precompact.
\\
\\
$(PS)_{*} \quad$If for some $M>0$, $\{u_j\}$ satisfies
\[
u_j\in E_j, \quad J(u_j)\leq M \quad \forall j, \quad ||(J|_{E_j})'(u_j)||_{E'_j} \to 0 \quad \text{as $j\to \infty$}
\]then $\{u_j\}$ is precompact.
These conditions follow from the fact that $J'$ is a compact perturbation of the Riesz representation map, and because such sequences are bounded. (Recall that in the proof of the Palais-Smale condition we only needed that $I_0(u_n)\leq C$. See equation (\ref{eq:7'})).
\\
\\
{\it Applying Tanaka's Theorem: The lower bound}\\
\\
The goal is to obtain a lower bound for $c_n$ that will contradict (\ref{eqApp29}). For $J(u)$ define the symmetric minimax levels
\begin{equation}\label{eq:39}
\beta_n := \inf_{g\in\Gamma_n}\sup_{u\in g(E_n)}J(u)
\end{equation}
Since $J(u)\leq I_0(u)$ by construction, we have $\beta_n \leq c_n$. So it will suffice to obtain a good lower bound on $\beta_n$. By \cite{Tanaka89} Theorem B, there exists a sequence $u_n$ such that

\begin{enumerate}
\item[i)] $J(u_n)\leq \beta_n$

\item[ii)] $J'(u_n) = 0$

\item[iii)] $n\leq index_0J''(u_n)$\\
\end{enumerate}
Where the extended Morse index $index_0J''(u)$ is the dimension of the maximal, negative semidefinite subspace corresponding to the form $J''(u)$. For simplicity we simply denote $u_n$ as $u$, holding $n$ fixed for the time being. Now
\[
J''(u)(v,w) = \langle v,w \rangle_{H^l_0(\Omega)} - \int_{\Omega}H''(u)vwdx
\]One basis for the maximal negative semidefinite subspace of this bilinear form is the set of eigenfunctions of $(-\Delta)^l-H''(u)$ with non-positive eigenvalues. So
\[
index_0J''(u) = \text{number of non-positive eigenvalues of $(-\Delta)^l-H''(u)$ on $L^2(\Omega)$}
\]By applying the Corollary of Proposition 1 we get
\begin{equation}\label{eq:40}
index_0J''(u) \leq C_{31}\int_{\Omega}\mathcal{B}(H''(u(x))dx+C_{31}
\end{equation}
So by Tanaka's theorem
\[
n \leq C_{32}\int_{\Omega}\mathcal{B}(H''(u(x))dx+C_{32}
\]where we take $n$ sufficiently large. Since $\Omega$ is of finite measure, the exact form of $H''(u)$ isn't important, only that it behaves like $(|u|+1)^{2\alpha_2-2}e^{(u^2+1)^b}$ for $|u|$ large. So that for some $C_{33}>0$ 
\[
\mathcal{B}(H''(u(x)) \leq C_{33}(|u|+1)^{3\alpha_2-2}e^{(u^2+1)^{b}}
\]
So
\[
\frac{n}{C_{34}} \leq \int_{\Omega}(|u|+1)^{3\alpha_2-2}e^{(u^2+1)^b}dx
\]Since $u=u_n$ is a critical point of $J$
\begin{eqnarray}\label{eq:41}
J(u)&=& A\int_{\Omega}\big[ bu^2(u^2+1)^{b-1}-1\big]e^{(u^2+1)^b}dx\nonumber\\
&\geq& C_{35}\int_{\Omega}(|u|+1)^{\alpha_2}e^{(u^2+1)^b}dx -C_{36} 
\end{eqnarray}
Let $\tau = (|u|+1)^{\alpha_2}e^{(u^2+1)^b}$, and for some $\gamma$ consider
\begin{eqnarray}\label{eq:42}
\kappa(\tau)&:=& \frac{\tau}{[\ln(\tau)]^{\gamma}}\nonumber\\
&\geq& C_{37}\big( |u|+1 \big)^{\alpha_2 -\alpha_2\gamma}e^{(u^2+1)^b}
\end{eqnarray}
We let $\alpha_2-\alpha_2\gamma = 3\alpha_2-2$, so $\gamma = 2/\alpha_2-2$. Note that $\gamma>0$ since $\alpha_2<1$. From (\ref{eq:42}) and the bound on $n$ we have
\begin{equation}\label{eq:43}
n \leq C_{39}\int_{\Omega}\kappa\bigg[ (|u|+1)^{\alpha_2}e^{(u^2+1)^b} \bigg]dx + C_{40}
\end{equation}
For large values of $\tau$, $\kappa''(\tau) <0$. So by the eventual concavity of $\kappa$ (\ref{eq:43}) gives
\begin{eqnarray}\label{eq:44}
n \leq C_{41}\kappa\bigg[\int_{\Omega} (|u|+1)^{\alpha_2}e^{(u^2+1)^b}dx \bigg] + C_{42}
\end{eqnarray}
Now apply inequality (\ref{eq:41}) by using the sublinearity of $\kappa$, inequality (\ref{eq:44}) gives
\begin{eqnarray}
n \leq C_{43}\cdot\kappa\big[ J(u_n)\big] + C_{44}\nonumber
\end{eqnarray} where we have included the subscript on $u$. Using the fact that $J(u_n)\leq b_n$, that $\kappa$ is eventually increasing, and that $b_n\to +\infty$ we have
\begin{eqnarray}\label{eq:45}
n \leq C_{43}\cdot\kappa[ b_n] + C_{44}
\end{eqnarray}
Let $\theta(\tau):= \tau[\ln(\tau)]^{\gamma}$, which is increasing and subexponential. Apply $\theta(\cdot)$ to both sides of (\ref{eq:45})
\begin{equation}\label{eq:46}
\theta(n) \leq C_{45}\cdot\theta \circ\kappa\big[ b_n\big] + C_{46}
\end{equation}
Now for large $\tau$
\[
\theta(\kappa(\tau)) =\tau\bigg[ 1- \gamma\frac{\ln\ln(\tau)}{\ln(\tau)} \bigg]^{\gamma} \leq \tau
\]So that for large $n$
\begin{equation}
C_{47}\theta(n)\leq \beta_n\nonumber
\end{equation}
i.e.
\begin{equation}
Cn[\ln(n)]^{\gamma}\leq \beta_n \leq c_n, \quad \gamma = 2/\alpha_2-2 \nonumber
\end{equation}
If $2/\alpha_2-2 > 1/\alpha_1$, as in the hypothesis of Theorem 1, this contradicts (\ref{eqApp29}) and so proves Theorem 1. $\blacksquare$

\section{Proof of the Eigenvalue bound}\label{sec5}

We essentially follow the argument in \cite{Solomyak94}. A general outline of the method is as follows. The initial step is to begin with a group of related embedding inequalities. These are then used to prove a theorem on piecewise-polynomial approximation. This approximation theorem then leads to an eigenvalue estimate for a compact operator on an appropriate Hilbert space. Finally one relates the eigenvalues of this operator to the non-positive eigenvalues of the Schr\"odinger operator of interest via the Birman-Schwinger principle.

Let $Q=(0,1)^d$ be the unit cube in $\mathbb{R}^d$ and let $u\in H^l(Q)=W^{l,2}(Q)$. Let $\Delta\subset\mathbb{R}^d$ be a parallelepiped with edges parallel to those of $Q$, and denote
\[
\mathcal{P}(l,d)= \text{vector space of all polynomials of degree $<l$ in $\Delta$,}
\]
\[
m(l,d):= \dim_{\mathbb{R}}\mathcal{P}(l,d).
\]
That is, we regard $\mathcal{P}(l,d)$ as that subspace of of $L^2(\mathbb{R}^d)$ consisting of functions supported in $\Delta$, and which in $\Delta$ are polynomials of degree less than $l$. We let $\mathcal{P}_{l,\Delta}$ be the corresponding orthonormal projection onto $\mathcal{P}(l,d)$. That is, $\mathcal{P}_{l,\Delta}$ is the $L^2$-orthogonal projection of $L^2(\mathbb{R}^d)$ onto $\mathcal{P}(l,d)$.

Furthermore, let $\Xi$ be a finite covering of Q by parallelepipeds $\Delta$. To any such covering and any $l>0$ we associate an operator of piecewise-polynomial approximation in $L^2(\mathbb{R}^d)$: For $\Xi = \{\Delta_j\}$, $1\leq j \leq card(\Xi)$, and with $\chi_j$ the characteristic function of the set $\Delta_j\setminus\cup_{i<j}\Delta_i$, we denote
\begin{equation}\label{eqSpec7}
K_{\Xi,l} = \sum_j\chi_j\mathcal{P}_{\Delta_j,l}.
\end{equation}
Note that $rank(K_{\Xi,l})\leq m(l,d)\cdot card(\Xi)$.

In this section we will need to recall the theory of Orlicz spaces (see \cite{Adams75},\cite{Kras-Rutic}, \cite{RaoRen}). Let $\mathcal{B},\mathcal{A}$ be a pair of mutually complementary $N$-functions, and $L_{\mathcal{B}}(\omega)$, $L_{\mathcal{A}}(\omega)$ be the corresponding Orlicz spaces on a set $\omega\subset\mathbb{R}^d$ of finite Lebesgue measure. We are primarily interested in the pair
\[
\mathcal{A}(t)=e^{|t|}-1-|t|, \quad \text{and}\quad \mathcal{B}(t)=(|t|+1)\ln(|t|+1)-|t|
\]
Then Solomyak's main theorem is
\begin{thm}[Solomyak \cite{Solomyak94}]\label{thm:4}
Let $Q=(0,1)^d$, $V\in L_{\mathcal{B}}(Q)$, $V\geq 0$. Then for any $n\in\mathbb{N}$ there exists a covering $\Xi=\Xi(V,n)$ of $Q$ by parallelepipeds $\Delta\subset Q$ such that
\begin{equation}\label{eqSpec8}
card(\Xi) \leq C_1n
\end{equation}
and for any $u \in H^l(Q)$, $2l=d$ we have
\begin{equation}\label{eqSpec9}
\int_QV|u-K_{\Xi,l}u|^2dx \leq C_2n^{-1}||V||_{\mathcal{B},Q}\int_Q|\nabla^lu|^2dx
\end{equation}
where $C_1, C_2$ depend only on $d$.

\end{thm}

Defined on $H^l_0(\Omega)$, we consider the quadratic form
\[
a(u,v) := \int_{\Omega}\nabla^lu\cdot \nabla^lvdx,
\]
with
\[
a(u):=a(u,u)= \int_{\Omega}|\nabla^lu|^2dx.
\]
As an unbounded quadratic form on $L^2(\Omega)$, $a(u)$ is symmetric and positive. It is also closed in $L^2(\Omega)$. To see this note that if $u_n\to u$ in $L^2(\Omega)$ and $a(u_m-u_n)\to 0$ then $\{u_n\}$ is Cauchy in $H^l_0(\Omega)$, and hence converges to some $\bar u$ in that space. By the generalized Poincar\'e Inequality $u_n\to \bar u$ in $L^2(\Omega)$. Thus $u=\bar u$ a.e. and so $u$ is (representable by) an element of $H^l_0(\Omega)$ and $\lim_{n\to\infty}a(u_n-u)=0$. On the space $L^2(\Omega)$ the unbounded operator $(-\Delta)^l$ is defined as the self-adjoint Friedrichs operator associated to the form $a(u)$. That is,

\begin{eqnarray}
D((-\Delta)^l) &:=& \{u\in H^l_0(\Omega): \text{the linear functional}\nonumber\\ &&v\mapsto \int \nabla^lu\cdot \nabla^lvdx \quad \text{is $L^2$-continuous,}\nonumber\\ &&\text{where}\quad v\in H^l_0(\Omega) \nonumber\}
\end{eqnarray}
and
\[
\langle (-\Delta)^lf,g \rangle_{L^2} = a(f,g)
\] for $f\in D((-\Delta)^l)$ and $g\in H^l_0(\Omega)$. We can do this since $a$ is closed. See for example section 5.5 in \cite{Weidmann}. Also, on $L^2(\Omega)$ we consider the form
\begin{equation}\label{eqSpec25}
a_V(u):= a(u) - \int_{\Omega}V(x)|u|^2dx
\end{equation} where $V\in L_{\mathcal{B}}$, $V\geq 0$. Here the form domain is $H^l_0(\Omega)\cap L^2(\Omega, Vdx)$. As we will see later in the proof of Theorem \ref{T}
\[
||u||_{L^2(\Omega, Vdx)}^2 \leq C_\Omega\cdot a(u) 
\]for $u\in H^l_0(\Omega)$. So the domain of $a_V$ is really just $d(a_V)= H^l_0(\Omega)$ and $H^l_0(\Omega)$ embeds into $L^2(\Omega, Vdx)$. As a matter of fact  for any $\epsilon>0$ there exists a constant $C(\epsilon)$ such that
\begin{equation}\label{eqSpec26}
||u||^2_{L^2(\Omega, Vdx)} \leq \epsilon a(u) + C(\epsilon)||u||^2_{L^2(\Omega)}.
\end{equation}
That is, the quadratic form 
\begin{equation}\label{eqSpec47}
\int_\Omega V(x)|u|^2dx
\end{equation}
has ``zero bound" relative to the form $a(u)$ in $L^2(\Omega)$. This is a consequence of Theorem \ref{T} below. By that theorem the operator corresponding to the quadratic form (\ref{eqSpec47}) in the Hilbert space $(H^l_0(\Omega),a(\cdot,\cdot))$ is compact. A well-known result then implies the zero-boundedness mentioned.

The bound (\ref{eqSpec26}) implies that $a_V(u)$ is lower semi-bounded in $L^2(\Omega)$. It is also closed in $L^2(\Omega)$. That is, for $\{u_n\}$ with $u_n\to u$ in $L^2(\Omega)$ and $a_V(u_n-u_m) \to 0$ then $a(u_n-u_m)\to 0$. So $\{u_n\}$ is Cauchy in $H^l_0(\Omega)$, and hence in $L^2(\Omega)$. Therefore $u\in H^l_0(\Omega)$. So as before we can define the associated self-adjoint Friedrichs operator on $L^2(\Omega, Vdx)$:
\[
A_V(u):= (-\Delta)^lu-\alpha V(x)u
\]whose domain is a subset of  $D((-\Delta)^l)$.

Suppose $A$ is a self-adjoint operator on a Hilbert space and that the spectrum of $A$ less than or equal to $\lambda\in\mathbb{R}$ is discrete. Then define $N(\lambda, A)$ to be the number of eigenvalues of $A$ less than or equal to $\lambda$, counted according to their multiplicity. For a compact, non-negative, symmetric operator $T$ denote by
\[
n(\lambda,T) = N(-\lambda, -T)
\]the number of eigenvalues of $T$ greater than or equal to $\lambda$. We now consider the quadratic functional
\[
b_V(u):= \int_{\Omega} V(x)|u(x)|^2dx
\]where $V\in L_{\mathcal{B}}$. If $b_V$ is bounded on $(H^l_0(\Omega), ||\cdot||_{H^l_0(\Omega)})$, then it generates a bounded, self-adjoint, non-negative operator on $H^l_0(\Omega)$ - say $T_V$. By definition
\begin{eqnarray}
u = T_Vf \iff && u\in H^l_0(\Omega);\nonumber\\
&& \int_{\Omega}\nabla^lu\cdot \nabla^lwdx = \int_{\Omega}Vfwdx,\qquad \forall w\in H^l_0(\Omega)\nonumber
\end{eqnarray}

\begin{thm}\label{T}
Let $\Omega \subset \mathbb{R}^{2l}$ be a bounded region with smooth boundary, and $V\in L_{\mathcal{B}}(\Omega)$. Then the operator $T_V$ is well-defined and compact on $H^l_0(\Omega)$, and there exists a constant $C_3=C_3(\Omega)$ such that for any $\lambda>0$
\begin{equation}\label{eq:51}
n(\lambda;T_V) \leq C_3||V||_{\mathcal{B},\Omega}\lambda^{-1}
\end{equation}
\end{thm}
{\bf Proof}: Let $Q\subseteq\mathbb{R}^{2l}$ be a cube such that $\bar\Omega\subseteq Q$. We can regard $Q$ as a unit cube, after rescaling. Let $W$ be the function on $Q$ equal to $V$ on $\Omega$ and $W=0$ on $Q\setminus \Omega$. Note that by the H\"older inequality for Orlicz spaces we write
\begin{eqnarray}\label{eq:52}
\int_{\Omega} V\cdot|u|^2dx \leq ||V||_{\mathcal{B},\Omega}\cdot ||u^2||_{\mathcal{A},\Omega}
\end{eqnarray} where the $N$-function $\mathcal{A}(t):=e^t-1-t$ is the Young function conjugate to $\mathcal{B}(t)$. Now by the well-known Sobolev-Orlicz embedding, see Theorem 8.25 in \cite{Adams75}, we have
\begin{proposition}\label{prop3}
There exists a constant $C$ such that for every $u\in H^l(\Omega)$, $2l=d$
\begin{equation}\label{eq:53}
||u^2||_{\mathcal{A},\Omega} \leq C ||u||_{W^{l,2}(\Omega)}^2
\end{equation}
\end{proposition}
Since $u \in H^l_0(\Omega)$ we can use the $H^l_0(\Omega)$ norm. So we get
\begin{equation}\label{eq:54}
\int_{\Omega} V\cdot|u|^2dx \leq C_d||V||_{\mathcal{B},\Omega}\cdot ||\nabla^lu||_{L^2}^2
\end{equation}
Thus $b_V$ is bounded as a quadratic form on $H^l_0(\Omega)$ (and on $\Pi\circ H^l_0(\Omega)\subseteq H^l_0(Q)$, where $\Pi$ is the natural extension operator). So $T_V$ is bounded on $H^l_0(\Omega)$ and 
\begin{equation}\label{eq:55}
n(\lambda;T_V)= 0 \quad \text{for $\lambda>C_d||V||_{\mathcal{B},\Omega}$}
\end{equation}
Now fix $\lambda\in(0,\lambda_0]$, where $\lambda_0 = C_2C_l||V||_{\mathcal{B},\Omega}$, for $C_2$ comes from Theorem \ref{thm:4}, and $C_l$ comes from
\[
||u||_{W^{l,2}}^2 \leq C_l||u||_{H^l_0}^2
\]for $u\in H^l_0(\Omega)$. Let $n$ be the minimal integer such that $n\lambda>\lambda_0$. For this $n$ and the function $W$, let $\Xi$ be the covering of $Q$ constructed in Theorem \ref{thm:4} and $K:= K_{\Xi,l}$ the corresponding operator (\ref{eqSpec7}).

For the subspace $\mathcal{F}:= ker(K\circ\Pi)$ of $H^l_0(\Omega)$
\[
codim \: \mathcal{F} \leq rank \: K \leq m(l,d)C_1\cdot n
\]
For $u\in H^l_0(\Omega)$ denote by $U:= \Pi(u)$. Then by Theorem (\ref{thm:4}) the following inequality holds
\begin{eqnarray}
\int_{\Omega}Vu^2dx &=& \int_QW|U-K( U)|^2dx\nonumber\\
&\leq&C_2n^{-1}||W||_{\mathcal{B},Q}\int_Q|D^lu|^2dx\nonumber\\
&\leq& C_lC_2n^{-1}||W||_{\mathcal{B},Q}||u||_{H^l_0}^2\nonumber\\
&\leq&\lambda ||u||_{H^l_0}^2\nonumber
\end{eqnarray}
This is enough to show that $T_V$ is compact. Indeed let $T_V|_{\mathcal{F}}$ denote the linear operator which is defined as $T_V$ on $\mathcal{F}$ and 0 on $\mathcal{F}^{\perp}$. Then $T_V = T_V|_{\mathcal{F}}+ T_V|_{\mathcal{F}^{\perp}}$ and $\mathcal{F}^{\perp}$ is finite dimensional. The above shows that $||T_V|_{\mathcal{F}}||\leq \lambda$. So taking $\lambda\to 0$, we see that $T_V$ is the limit of $T_V|_{\mathcal{F}^{\perp}}$ in the uniform norm. Therefore $T_V$ is compact and its spectrum consists of eigenvalues. If $\{u_j\}$ is an eigenvector with eigenvalue $\lambda_j \geq \lambda$ then

\begin{eqnarray}\label{eq:56}
\lambda\leq \lambda_j = \frac{\langle T_Vu_j,u_j \rangle_{H^l_0}}{\langle u_j,u_j \rangle_{H^l_0}}
=\frac{\int_\Omega V|u_j|^2dx}{\int_\Omega |\nabla^lu_j|^2dx}
\end{eqnarray}
So $u_j\notin \mathcal{F}$. Since eigenvectors are orthogonal
\begin{eqnarray}\label{eq:57}
n(\lambda;T_V) \leq codim\: \mathcal{F} &\leq& m(l,d)C_1n\nonumber\\
&<& m(l,d)C_1\bigg(\frac{\lambda_0}{\lambda}+1\bigg)\nonumber\\
&\leq& 2m(l,d)C_1\frac{\lambda_0}{\lambda}, \quad \lambda \leq \lambda_0
\end{eqnarray}
The required estimate (\ref{eq:51}), with $C_3 = 2m(l,d)C_1\max\{C_lC_2,C_d\}$, where $C_d$ is given in (\ref{eq:54}) and $C_l$ is given in the inequality after (\ref{eq:55}), is a consequence of (\ref{eq:55}) and (\ref{eq:57}):
\\
i) If $C_lC_2 \geq C_d$ then (\ref{eq:57}) gives the result.
\\
ii) If $\lambda\leq \lambda_0 = C_lC_2||V||_{\mathcal{B},\Omega}$, again (\ref{eq:57}) gives the result.
\\
iii) If $C_lC_2||V||_{\mathcal{B},\Omega} < \lambda \leq C_d||V||_{\mathcal{B},\Omega}$ then (\ref{eq:57}) is applied to $\tilde{\lambda}= \lambda\frac{C_lC_2}{C_d}$ and that gives the result.
\\
$\blacksquare$
\\
\\
{\bf Proof of Proposition 1 : The Briman-Schwinger Principle}
\\
\\
The reference here is section 1 in \cite{BirmanSolomyak3}. Let $a(u)$ be a positive, symmetric, and closed quadratic form in a Hilbert space $\mathcal{H}$ with domain $D(a)\subset \mathcal{H}$. Let $b(u)$ be another non-negative, symmetric quadratic form such that
\begin{equation}\label{eq:58}
b(u)\leq C\cdot a(u), \qquad u\in d
\end{equation}
Consider the space $\tilde{D}(a)$ - the completion of $D(a)$ in the inner product given by $a(\cdot,\cdot)$. By (\ref{eq:58}) $b$ can be extended to all of $\tilde{D}(a)$. The extended form defines on $\tilde{D}$ a bounded, self-adjoint, non-negative operator, which we denote by $B:(\tilde{D}(a),a(\cdot,\cdot))\rightarrow  (\tilde{D}(a),a(\cdot,\cdot))$.

\begin{proposition}[Birman-Schwinger Principle]
Suppose (\ref{eq:58}) is satisfied and the operator $B$ is compact as an operator from $(\tilde{D}(a),a(\cdot,\cdot))$ to itself. Then for any $\alpha>0$ the quadratic form
\[
a_{\alpha}(u):= a(u)-\alpha b(u), \qquad u\in D
\]is semi-bounded from below and closed in $\mathcal{H}$. As usual, this implies that there is a corresponding self-adjoint Friedrichs operator $A_{\alpha b}$ associated with this form. For $A_{\alpha b}$ the non-positive spectrum is finite and
\begin{equation}\label{eq:59}
N(0;A_{\alpha b}) = n(\alpha^{-1};B)
\end{equation}
\end{proposition}
This result comes from the variational characterization of $N(0; A)$, often referred to as Glazman's lemma. For us $\mathcal{H} = L^2(\Omega)$, $D(a)=\tilde{D}(a)=H^l_0(\Omega)$, $a(u)=\int_{\Omega}|\nabla^lu|^2dx$, $b(u)=\int_{\Omega}V|u|^2dx$, $\alpha=1$, $B=T_V$. From the proof of Theorem 4 (\ref{eq:58}) is satisfied, $T_V$ is compact, and the form $a(u)$ is closed in $L^2(\Omega)$. We thus obtain
\[
N(0,(-\Delta)^l-V(x)) = n(1, T_V)
\] So by Theorem 4
\[
N(0,(-\Delta)^l-V(x)) \leq C_3||V||_{\mathcal{B},\Omega}
\] 
$\blacksquare$

\section{Problem (R), the radial problem on an annulus}\label{sec6}

{\it 6.1 The problem and its variational setup}
\\

Here $\Omega = A_{R_0}^R := \{x\in \mathbb{R}^{2l}: R_0 < |x| < R \}$ will denote an annulus, with $R_0>0$ and $R<+\infty$. We seek radial solutions to the problem 

\[
\tag{R}
\begin{cases}
(-\Delta)^lu = 2ue^{u^2} + \varphi(x,u) \quad \mbox{in $\Omega$} 
\\
\\
\bigg(\frac{\partial}{\partial \nu}\bigg)^ju\Bigg|_{\partial\Omega} = 0, \quad j = 0, \dots, l-1
\end{cases}
\]
where $\varphi(x,u)=\varphi(|x|,u)$. The proper space for this problem is
\[
H_r := \{ u \in H_0^l(\Omega): u(x)=u(|x|) \ \text{a.e. in $\Omega$}\}
\] and the corresponding functional is
\[
I_1(u):= \frac{1}{2}\int_{\Omega}|\nabla^lu|^2dx - \int_{\Omega} (e^{u^2}-1)dx - \int_{\Omega}\Phi(x,u) dx
\] where $\Phi(x,u)=\int_{0}^{u}\varphi(x,t)dt$ is not even in $u$. Since $\varphi(x,u)$ is radial in its explicit dependence on $x$, critical points of $I_1$, even when restricted to $H_r$, still correspond to generalized solutions of (R). This can be seen from a simple direct calculation using spherical coordinates. The general principle behind this fact is called the Principle of Symmetric Criticallity, see \cite{Palais79}. Concerning the size of the perturbation, we assume there exists  $\beta<1$ and $C>0$ such that
\begin{equation}\label{eq:0060}
|\Phi(x,t)| + |\varphi(x,t)t| \leq Ct^2e^{\beta t^2}
\end{equation} for all $x\in\Omega$, and $t\in\mathbb{R}$ with $|t|$ large. Similarily to before, the path of functionals we will be concerned with is
\[
I_{\theta}(u):= \frac{1}{2}\int_{\Omega}|\nabla^lu|^2dx - \int_{\Omega} (e^{u^2}-1)dx - \theta\int_{\Omega}\Phi(x,u) dx
\]where $\theta \in [0,1]$.
\\
\\
{\it 6.2 Bolle's and Tanaka's requirements}
\\
\\
The compactness of Palais-Smale sequences in the space $H_r$ as required by Bolle's condition (H1) follows as in the general case, with the exception that instead of Adam's inequality we use Lemma \ref{lemma3}, given below in section \ref{sec7}. Thus the radial setting allows us to consider an nonlinearity which was of critical growth in the unrestricted setting earlier. For the condition (H2) assume $|I_{\theta}(u)|\leq b$. In fact, we only need to assume that $I_\theta(u)\leq b$:
\begin{eqnarray}
||I_{\theta}'(u)||_{H^{-l}}||u||_{H^l} \geq -\langle I_{\theta}'(u),u \rangle = -\int|\nabla^lu|^2dx + \int 2u^2e^{u^2}+\theta\varphi(x,u)udx \nonumber\\
= \int|\nabla^lu|^2dx + \int (2u^2-4)e^{u^2}+\theta\varphi(x,u)u-4\theta\Phi(x,u)dx - 4I_{\theta}(u)\nonumber - C_1\\
\geq ||u||^2_{H^l(\Omega)} + c\int u^2e^{u^2}dx - C_2 - b\nonumber\\
\geq c\bigg|\int \Phi(x,u)dx\bigg| -C_3= c\bigg|\frac{\partial}{\partial\theta}I_{\theta}(u)\bigg|-C_3\nonumber
\end{eqnarray}
where $c$ is some small positive constant, and where we used (\ref{eq:0060}) in the last inequaliy. This verifies condition (H2). Next for condition (H3), let $u$ be a critical point of $I_{\theta}$. Then
\begin{eqnarray}
I_{\theta}(u) =I_\theta(u)-\frac{1}{2}\langle I'_\theta(u),u\rangle\nonumber\\
= \int (u^2-1)e^{u^2} + \frac{\theta}{2}\varphi(x,u)u-\theta\Phi(x,u)dx-C_4\nonumber\\
\geq c\int u^2e^{u^2}dx -C_5\nonumber
\end{eqnarray}where again $c$ is some small positive constant. Applying Jensen's inequality and (\ref{eq:0060}) we have for sufficiently large constants $C_6$ and $C_7$
\begin{eqnarray}
C_6[|I_{\theta}(u)|+1]^{\beta} \geq \int\bigg( u^2e^{u^2} \bigg)^{\beta}dx+C_7\nonumber\\
\geq \bigg|\int \Phi(x,u)dx\bigg| = \bigg| \frac{\partial}{\partial\theta}I_{\theta}(u) \bigg|
\end{eqnarray}
So condition (H3) holds with $\bar{f}_1(\theta,t)=\bar{f}_2(\theta,t) = f(t):=C[|t|+1]^{\beta}$. The condition (H4) follows as before. If we assume that only the second possibility of Theorem \ref{thm:3} holds for sufficiently large $n\in\mathbb{N}$ then for some $K>0$ such that
\[
c_{n+1}-c_n \leq K(\bar{f}_1(c_{n+1})+\bar{f}_2(c_n)+1)
\]for sufficiently large $n$. More consicely, by enlarging $K$ if necessary, this means
\[
c_{n+1}-c_n \leq K((c_{n+1})^\beta+(c_n)^\beta+1)
\]for sufficiently large $n$. Finally, using the fact that $\beta<1$, this implies that for some $A>0$
\[
c_n \leq An^{\frac{1}{1-\beta}} \quad \text{for $n>n_0$}
\] with $n_0$ sufficiently large. The argument is the same as that used for (\ref{eqApp29}).
\\
\\
\\
{\it 6.2 The lower bound}
\\
\\
As before, the goal is to obtain a lower bound for $c_n$ that will contradict the alleged upper bound. The requirements in Tanaka's Theorem are all verified as before, with no new phenomena appearing. By Tanaka's theorem there exists a sequence $u_n$ in the Hilbert space $H_r$ such that

\begin{enumerate}
\item[i)] $I_0(u_n)\leq c_n$
\item[ii)] $I_0'(u_n)=0$
\item[iii)] $n \leq index_0I_0''(u_n)$
\end{enumerate}
For simplicity we denote $u_n$ as $u$, holding $n$ fixed for the moment. As before,  
\[
index_0I_0''(u) = \text{number of non-positive eigenvalues of $(-\Delta)^l - 2(u^2+2)e^{u^2}$}
\]By applying Proposition \ref{prop:2} we get

\begin{eqnarray}
n^{2l} \leq [index_0I_0''(u)]^{2l} &\leq& C_8 \int_{A_{R_0}^R} 2(u^2+2)e^{u^2}\log \bigg(\frac{R}{|x|}\bigg)^2idx\nonumber \\
& \leq& C_{R_0} \int_{A_{R_0}^R} 2(u^2+2)e^{u^2}dx\nonumber
\end{eqnarray} 
Since $u=u_n$ is a critical point of $I_0$ as before we have that $c_n \geq I_0(u_n) \geq c\int u^2e^{u^2}dx-C_9$. Therefore the we obtain
\[
c_n \geq C_{10}\cdot n^{2l}
\] for sufficiently large n. Therefore if $2l> \frac{1}{1-\beta}$ this contradicts the upper bound, and proves Theorem \ref{thm:02}.

\section{Piecewise-polynomial approximation and spectral estimates in the annular case}\label{sec7}

Here we prove Proposition \ref{prop:2}. As mentioned earlier, we first need some appropriate inequalities in the radial case to take the place of the Orlicz-Sobolev inequality of Proposition \ref{prop3}. Lemma \ref{lemma3} below is a generalization of an inequality by Ni, see \cite{Ni1}.

\begin{lemma}\label{lemma3}
Let $u \in H_r$.
\begin{enumerate}
\item[a)] If $d=2l=2$ then
	\[
	|u(x)| \leq C_{d}||\partial_ru||_{L^2} \cdot \sqrt{\log\bigg( \frac{R}{|x|} \bigg)}
	\]for $x\in \Omega = A_{R_0}^R$.
\item[b)] If $d=2l > 2$ then
\[
	|u(x)| \leq C_{d}||\partial_r^lu||_{L^2} \cdot \log\bigg( \frac{R}{|x|} \bigg)
	\] for $x\in\Omega = A_{R_0}^R$.
\end{enumerate}
\end{lemma}

{\bf Proof}:\\
\\
a) For simplicity we write $u = u(r)$ as a function of the radial variable. By a density argument we may assume that $u\in C^{\infty}_0(\Omega)$. For $r\in [R_0,R]$

\begin{equation*}
-u(r) = u(R)-u(r) = \int^{R}_{r}u(\rho)d\rho
\end{equation*}so

\begin{eqnarray}
|u(r)| &\leq& \int^{R}_{r}|u(\rho)|d\rho \nonumber\\
&\leq& \bigg( \int^{R}_{r}|u(\rho)|^2\rho d\rho \bigg)^{1/2}\bigg( \int^{R}_r \frac{1}{\rho} d\rho \bigg)^{1/2} \nonumber\\
&\leq& C_{d}||\partial_ru||_{L^2} \cdot \log\bigg( \frac{R}{r} \bigg)^{1/2}\nonumber
\end{eqnarray}
which is the required result.
\\
\\
b) Again we take $u \in C^{\infty}_0$.

\begin{eqnarray}
u(r) &=& u(r)-u(R) = -\int_{r}^{R}u'(\rho_1)d\rho_1 \nonumber\\
	&=& \int_{r}^{R}u'(R)-u'(\rho_1)d\rho_1 = \int_{r}^{R}\int_{\rho_1}^{R}u''(\rho_2)d\rho_2d\rho_1 \nonumber\\
	&\vdots& \nonumber\\
	&=& (-1)^l\int_{r}^{R}\int_{\rho_1}^{R} \cdots \int_{\rho_{l-1}}^{R}u^{(l)}(\rho_l)d\rho_l \cdots d\rho_1\nonumber
\end{eqnarray}So

\begin{eqnarray}
|u(r)| &\leq& \int_{r}^{R}\int_{\rho_1}^{R} \cdots \int_{\rho_{l-1}}^{R}|u^{(l)}(\rho_l)|d\rho_l \cdots d\rho_1 \nonumber\\
	&=& \int_{r}^{R}\int_{\rho_1}^{R} \cdots \int_{\rho_{l-1}}^{R}|u^{(l)}(\rho_l)|\rho_l^{\frac{d-1}{2}}\rho_l^{\frac{1-d}{2}}d\rho_l \cdots d\rho_1 \nonumber\\
	&\leq& \int_{r}^{R}\int_{\rho_1}^{R} \cdots \int_{\rho_{l-2}}^{R} \bigg( \int_{\rho_{l-1}}^{R}|u^{(l)}(\rho_l)|^2\rho_l^{d-1}d\rho_l  \bigg)^{1/2}  \bigg( \int_{\rho_{l-1}}^{R} \rho_{l}^{1-d}d\rho_l  \bigg)^{1/2}d\rho_{l-1} \cdots d\rho_1\nonumber
\end{eqnarray}

\begin{eqnarray}
	&\leq& C_d ||\partial_r^lu||_{L^2} \cdot \int_{r}^{R}\int_{\rho_1}^{R} \cdots \int_{\rho_{l-2}}^{R}\rho_{l-1}^{1-l}d\rho_{l-1} \cdots d\rho_1\nonumber\\
	&\leq& C_d ||\partial_r^lu||_{L^2} \cdot \int_{r}^{R}\frac{1}{\rho_1}d\rho_1\nonumber\\
	&=&C_d ||\partial_r^lu||_{L^2} \cdot \log\bigg( \frac{R}{r} \bigg)\nonumber
\end{eqnarray} which is the required estimate. $\blacksquare$\\
\\

When proving a radial version of Theorem \ref{thm:4} it is necessary to have at one's disposal inequalities of the above type, but without the zero boundary conditions. The key is to find the appropriate $(l-1)$th-degree polynomial to subtract from $u$, so that the remainder can be controlled by the $l$th-order derivative of $u$. It is not surprising that this is the same polynomial approximation which appears in the radial version of Theorem \ref{thm:4}.

\begin{lemma}\label{lemma4}
Let $u \in H^l(A)$, where $A=A_{R_0}^R$ is an annulus centered at the origin in $\mathbb{R}^{2l}$, and $u(x) = u(|x|)$ a.e in $A$.
\begin{enumerate}
\item[a)] If $d=2l=2$ then
	\[
	|u(x)-\bar{u}_A| \leq C_{d}||\partial_ru||_{L^2(A)} \cdot \bigg[1+\log\bigg( \frac{R}{|x|} \bigg)^{1/2}\bigg]
	\]for $x\in A$, where $\bar{u}_A$ is the average value of $u$ in $A$. 
\item[b)] When $d=2l > 2$. First define
	\[
	\tau_l(u)(r,s) := \sum_{n=0}^{l-1}\frac{u^{(n)}(s)}{n!}(r-s)^n
	\]and
	
	\[
	P_{l,A}(u)(r) := \frac{1}{|A|}\int_A \tau_l(u)(r,|x|)dx
	\]which is a polynomial in $r$ of degree $\leq l-1$ and linear in $u$. Then
	
	\[
	|u(x)-P_{l,A}(u)(|x|)| \leq C_{d}||\partial_r^lu||_{L^2(A)} \cdot \bigg[1+\log\bigg( \frac{R}{|x|} \bigg) \bigg]
	\] for $x\in A=A_{R_0}^R$.
\end{enumerate}
\end{lemma}

{\bf Proof}: We assume $u\in C^{\infty}(A)$ and radially symmetric. Then
\begin{equation}\label{eq:60}
|u(r)-\bar{u}_A| = \bigg| \frac{1}{|A|}\int_A u(r)-u(x)dx \bigg|\leq \frac{1}{|A|}\int_A |u(r)-u(x)|dx
\end{equation} Now
\[
u(r)-u(x) = \int_{|x|}^{r}u'(\rho)d\rho
\]So 
\[
|u(r)-u(x)| \leq \int_{|x|,r}|u'(\rho)|\rho^{1/2}\rho^{-1/2}d\rho
\]where the notation $\int_{a,b}$ denotes unoriented integration over the interval with endpoints $a$ and $b$. So
\begin{eqnarray*}
|u(r)-u(x)| &\leq& \bigg( \int_{|x|,r}|u'(\rho)|^2\rho d\rho \bigg)^{1/2}\bigg( \int_{|x|,r}\rho^{-1}d\rho  \bigg)^{1/2}\\
	&\leq& C_d ||\partial_ru||_{L^2(A)} \cdot \bigg|\log\bigg(\frac{r}{|x|} \bigg) \bigg|^{1/2}
\end{eqnarray*}
Plugging this into the earlier inequality gives
\begin{eqnarray*}
|u(r)-\bar{u}_A| &\leq& \frac{1}{|A|}\int_A C_d ||\partial_ru||_{L^2(A)} \cdot \bigg|\log\bigg(\frac{r}{|x|} \bigg) \bigg|^{1/2}dx \\
	&=&\frac{C_d}{|A|}||\partial_ru||_{L^2(A)}\cdot \int_{R_0}^R  \bigg|\log\bigg(\frac{r}{\rho} \bigg) \bigg|^{1/2}\rho d\rho
\end{eqnarray*}
It's easy to check that the value of the integral is bounded above by a constant times
\[
R(R-R_0)\bigg[1+ \log\bigg( \frac{R}{r} \bigg)^{1/2} \bigg].
\]
To see this we evaluate the integral in two parts:
\[
I_1 = \int_{R_0}^{r}\log\bigg(\frac{r}{\rho}\bigg)^{1/2}\rho d\rho
\]
and
\[
I_2 = \int_{r}^{R}\log\bigg(\frac{\rho}{r}\bigg)^{1/2}\rho d\rho.
\]
In $I_1$ we let $t=\log(r/\rho)$ and so
\[
I_1 = r^2\int_{0}^{\log(r/R_0)}t^{1/2}e^{-2t}dt \leq C_0r^2\int_{0}^{\log(r/R_0)}e^{-t}dt = C_0r(r-R_0)
\]
\[
\leq C_0R(R-R_0).
\]
For $I_2$ we simply notice 
\[
I_2 = \int_{r}^{R}\log\bigg(\frac{\rho}{r}\bigg)^{1/2}\rho d\rho \leq \log\bigg(\frac{R}{r}\bigg)^{1/2}R(R-R_0).
\]
Since $A$ is a $d=2$ dimensional annulus we have that
\[
\frac{R(R-R_0)}{|A|} \leq C_d
\]
Thus
\[
|u(r)-\bar{u}_A| \leq C_d ||\partial_ru||_{L^2(A)} \cdot \bigg[1+ \log\bigg( \frac{R}{r} \bigg)^{1/2} \bigg]
\] which is the required result.
\\
\\
b) As before, we assume $u$ is smooth. Then 
\begin{eqnarray}\label{eq:(61)}
|u(x)-P_{l,A}(u)(|x|)| &=& \bigg | \frac{1}{|A|}\int_A u(r)-\tau_l(u)(r,|x|)dx \bigg| \nonumber \\
	&\leq& \frac{1}{|A|}\int_A |u(r)-\tau_l(u)(r,|x|)|dx
\end{eqnarray}
Set $v(r):= u(r)-\tau_l(u)(r,|x|)$. Note that, when keeping $|x|$ fixed, we have
\[
v^{(n)}(r)|_{r=|x|} = 0 \quad \text{for $0 \leq n \leq l-1$}
\] where the differentiation is partial differentiation w.r.t. $r$. We seek to estimate $v(r)$ by repeatedly applying this property.
\begin{eqnarray*}
v(r) &=& v(r)-v(|x|) = \int_{|x|}^{r}v'(\rho_1)d\rho_1 \\
	&=&\int_{|x|}^{r}v'(\rho_1)-v'(|x|)d\rho_1 = \int_{|x|}^{r}\int_{|x|}^{\rho_1}v''(\rho_2)d\rho_2d\rho_1 \\
	&=& \int_{|x|}^{r}\int_{|x|}^{\rho_1}v''(\rho_2)-v''(|x|)d\rho_2d\rho_1 \\
	&\vdots& \\
	&=&\int_{|x|}^{r}\int_{|x|}^{\rho_1} \cdots \int_{|x|}^{\rho_{l-1}}  v^{(l)}(\rho_l)d\rho_l \cdots d\rho_1
\end{eqnarray*}
thus

\[
|v(r)| \leq \int_{|x|,r}\int_{|x|,\rho_1} \cdots \int_{|x|,\rho_{l-1}}  |v^{(l)}(\rho_l)|d\rho_l \cdots d\rho_1
\]
We apply H\"older's inequality to the inner most integral
\begin{eqnarray*}
|v(r)| \leq \int_{|x|,r}\int_{|x|,\rho_1} \cdots \int_{|x|,\rho_{l-2}} \bigg(  \int_{|x|,\rho_{l-1}} |v^{(l)}(\rho_l)|^2\rho_l^{d-1} d\rho_l \bigg)^{1/2}\times \\ \bigg( \int_{|x|,\rho_{l-1}} \rho_l^{1-d}d\rho_l \bigg)^{1/2}d\rho_{l-1} \cdots d\rho_1
\end{eqnarray*}

\[
\leq C_d ||\partial_r^{l}u||_{L^2(A)}\cdot \int_{|x|,r}\int_{|x|,\rho_1} \cdots \int_{|x|,\rho_{l-2}} |\rho_{l-1}^{2-d}-|x|^{2-d}|^{1/2}d\rho_{l-1} \cdots d\rho_1
\]
after using $v^{(l)}(\rho_l)= \partial_r^lu(\rho_l)$. Returning to \ref{eq:(61)} we get
\begin{eqnarray}\label{eq:62}
|u(r)-P_{l,A}(u)(r)|\leq \frac{C_d}{|A|}||\partial_r^{l}u||_{L^2(A)}\times&& \nonumber \\ \int_A \int_{|x|,r}\int_{|x|,\rho_1} \cdots \int_{|x|,\rho_{l-2}} |\rho_{l-1}^{2-d}-|x|^{2-d}|^{1/2} && d\rho_{l-1} \cdots d\rho_1 dx
\end{eqnarray}
We seek to estimate the above integral. First, by converting to polar coordinates the integral becomes (after factoring out a constant depending only on the dimension)
\[
\int_{R_0}^{R}\int_{\rho_0,r}\int_{\rho_0,\rho_1} \cdots \int_{\rho_0,\rho_{l-2}} |\rho_{l-1}^{2-d}-\rho_0^{2-d}|^{1/2} d\rho_{l-1} \cdots d\rho_1 \rho_0^{d-1}d\rho_0
\]
We divide the integration by $d\rho_0$ into two pieces. One where $R_0 \leq \rho_0 \leq r$ and one where $r\leq \rho_0 \leq R$. This allows us to properly orient the endpoints. The first integral is
\[
\int_{R_0}^{r}\int_{\rho_0}^{r}\int_{\rho_0}^{\rho_1} \cdots \int_{\rho_0}^{\rho_{l-2}}\big( \rho_0^{2-d} - \rho_{l-1}^{2-d} \big)^{1/2} d\rho_{l-1} \cdots d\rho_1 \rho_0^{d-1} d\rho_0
\]
\[
\leq \int_{R_0}^{r}\int_{\rho_0}^{r}\int_{\rho_0}^{\rho_1} \cdots \int_{\rho_0}^{\rho_{l-2}}\rho_0^{1-d/2} d\rho_{l-1} \cdots d\rho_1 \rho_0^{d-1} d\rho_0
\]
\[
= \int_{R_0}^{r}\int_{\rho_0}^{r}\int_{\rho_0}^{\rho_1} \cdots \int_{\rho_0}^{\rho_{l-3}}\rho_0^{1-d/2} (\rho_{l-2}-\rho_0) d\rho_{l-2} \cdots d\rho_1 \rho_0^{d-1} d\rho_0
\]
\[
\leq \int_{R_0}^{r}\int_{\rho_0}^{r}\int_{\rho_0}^{\rho_1} \cdots \int_{\rho_0}^{\rho_{l-3}}\rho_0^{1-d/2} \rho_{l-2} d\rho_{l-2} \cdots d\rho_1 \rho_0^{d-1} d\rho_0
\]
\[
\vdots
\]
\[
\leq C_l \int_{R_0}^{r}\int_{\rho_0}^{r}\rho_0^{1-d/2}\rho_1^{l-2}d\rho_1 \rho_0^{d-1} d\rho_0
\]
\[
\leq C\int_{R_0}^{r}\rho_0^{1-d/2}r^{l-1}\rho_0^{d-1} d\rho_0
\]
\[
= C\int_{R_0}^{r}\rho_0^{d/2}r^{l-1} d\rho_0
\]
\[
\leq Cr^{d/2 +l -1}(r-R_0)
\]
\begin{equation}\label{eq:63}
\leq CR^{d/2 +l -1}(R-R_0) = CR^{d-1}(R-R_0)
\end{equation}where we have used $d = 2l$ in the final equality. \\
\\
The second integral is 

\[
\int_{r}^{R}\int_{r}^{\rho_0}\int_{\rho_1}^{\rho_0} \cdots \int_{\rho_{l-2}}^{\rho_0}\big(  \rho_{l-1}^{2-d} - \rho_0^{2-d} \big)^{1/2} d\rho_{l-1} \cdots d\rho_1 \rho_0^{d-1}d\rho_0
\]
\[
\leq \int_{r}^{R}\int_{r}^{\rho_0}\int_{\rho_1}^{\rho_0} \cdots \int_{\rho_{l-2}}^{\rho_0}\rho_{l-1}^{1-d/2} d\rho_{l-1} \cdots d\rho_1 \rho_0^{d-1}d\rho_0
\]
\[
\vdots
\]
\[
\leq C_d\int_{r}^{R}\int_{r}^{\rho_0}\rho_1^{l-1-d/2}d\rho_1 \rho_0^{d-1}d\rho_0
\]
Now using $l-1-d/2=-1$,
\begin{equation}\label{eq:64}
= C_d \int_{r}^{R}\log\big( \frac{\rho_0}{r} \big) \rho_0^{d-1} d\rho_0 \leq C_d R^{d-1}(R-R_0)\log\bigg( \frac{R}{r} \bigg)
\end{equation}
Combining \ref{eq:62}, \ref{eq:63}, and \ref{eq:64} gives
\begin{equation}\label{eq:65}
|u(r)-P_{l,A}(u)(r)|\leq C_d \frac{R^{d-1}(R-R_0)}{|A|}||\partial_r^lu||_{L^2(A)}\bigg[ 1+ \log \bigg( \frac{R}{r} \bigg)\bigg]
\end{equation}
Since $A$ is an annulus, $\frac{R^{d-1}(R-R_0)}{|A|} \leq C_d$. So
\begin{equation}\label{eq:66}
|u(r)-P_{l,A}(u)(r)|\leq C_d ||\partial_r^lu||_{L^2(A)}\bigg[ 1+ \log \bigg( \frac{R}{r} \bigg)\bigg]
\end{equation}
Which is the required estimate. $\blacksquare$\\
\\

We will apply Lemma \ref{lemma4} to our main region $\Omega = A_{R_0}^R$, where $R<\infty$ and $R_0>0$ are, respectively, the fixed outer and inner radii of $\Omega$. Let $\tilde{u}(x)$, $x\in\Omega$, be as in that lemma. We have that for $\tilde{r}\in [R_0,R]$
\[
|\tilde{u}(\tilde{r})-P_{l,\Omega}\tilde{u}(\tilde{r})|^2 \leq c_d \bigg[ 1+ \log\bigg( \frac{R}{\tilde{r}} \bigg) \bigg]^{2i}\int_{R_0}^{R}|\partial_{\tilde{r}}^l\tilde{u}|^2\tilde{\rho}^{d-1}d\tilde{\rho}
\]
So
\begin{equation}\label{eq:67}
|\tilde{u}(\tilde{r})-P_{l,\Omega}\tilde{u}(\tilde{r})|^2 \leq \kappa_0\int_{R_0}^{R}|\partial_{\tilde{r}}^l\tilde{u}|^2\tilde{\rho}^{d-1}d\tilde{\rho}
\end{equation}
where
\[
\kappa_0 = \kappa_0(d,R_0,R)= c_d \bigg[ 1+ \log\bigg( \frac{R}{R_0} \bigg) \bigg]^{2i}
\]
Now consider a change of variable, replacing the domain $\Omega=A_{R_0}^R$ with a smaller annulus contained in it, $A$. It is centered at the origin, with inner radius $R_A$, and outer radius $R^A$, $R_0\leq R_A < R^A \leq R$:\\
\\
Let $r = \mu \tilde{r} +\beta$, where $\mu := \frac{R^A-R_A}{R-R_0}$ and $\beta := \frac{R_0R^A-RR_A}{R^A-R_A}$. For a radially symmetric $\tilde{u}\in H^l(\Omega)$, let $u(r):=\tilde{u}(\tilde{r})$. Define 

\begin{equation}\label{eq:068}
P_{l,A}u(r) := P_{l,\Omega}\tilde{u}(\tilde{r})
\end{equation}
Inequality (\ref{eq:67}) becomes
\begin{equation}\label{eq:069}
|u(r)-P_{l,A}u(r)|^2 \leq \kappa_0 \mu^{2l-1}\int_{R_A}^{R^A}|\partial_r^lu|^2\tilde{\rho}d\rho
\end{equation} where $\rho = \mu\tilde{\rho}+\beta$. Since $R_0>0$
\[
\frac{R_0}{R} \leq \frac{\tilde{\rho}}{\rho} \leq \frac{R}{R_0}
\]
This and inequality (\ref{eq:069}) give
\begin{equation}\label{eq:070}
|u(r)-P_{l,A}u(r)|^2 \leq \kappa' \mu^{2l-1}\int_A|\partial_r^lu|^2dx, \quad \kappa' := \kappa_0\big(R/R_0\big)^{d-1} 
\end{equation}
But clearly
\[
\mu^{2l-1} \leq \bigg( \frac{c_d}{R_0^{d-1}(R-R_0)} \bigg)^{2l-1}|A|^{2l-1}
\]So finally, inequality (\ref{eq:070}) gives
\begin{equation}\label{eq:071}
|u(r)-P_{l,A}u(r)|^2 \leq \kappa'' |A|^{2l-1}\int_A|\partial_r^lu|^2dx, \quad \kappa'' := \kappa'\bigg( \frac{c_d}{R_0^{d-1}(R-R_0)} \bigg)^{2l-1}
\end{equation} for all radially symmetric $u\in H^l(A)$. We will use this result as the basis for establishing a radial analogue of Theorem \ref{thm:4}. Before we proceed with that, we need a lemma on functions of sets.

Let $\mathcal{J}$ be a nonnegative function of half-open annuli $A\subseteq \Omega$ (always taken to be centered at the origin), which is super-additive. That is, if an annulus $A$ is partitioned into finitely many annuli $\{A_j\}$, then $\sum_j\mathcal{J}(A_j)\leq \mathcal{J}(A)$. The $\mathcal{J}$ that we are interested in is
\[
\mathcal{J}(A):= \int_AV(x)\bigg[1+ \log\bigg(\frac{R}{|x|}\bigg) \bigg]^{2i}dx
\]
For a partition $\Xi=\{A_j\}$ of $\Omega$ into (half-open) annuli define
\[
G(\mathcal{J},\Xi):= \max_{A\in\Xi}|A|^{2l-1}\mathcal{J}(A)
\]Then by Theorem 1.5 in \cite{BirmanSolomyak2} we have
\begin{lemma}\label{lemma5}
For any natural number $n$ there exists a partition $\Xi$ of $\Omega$ into (half-open) annuli such that
\[
card(\Xi) \leq n \quad \text{and}
\]
\[
G(\mathcal{J},\Xi) \leq C(l,d)n^{-2l}\mathcal{J}(\Omega).
\]
\end{lemma}Its proof is given in section 2.2 of that reference. Let $\Xi$ be the partition guaranteed by the lemma. For a set $A$ let $\chi_A$ denote its characteristic function. Define
\[
K_{\Xi,l} := \sum_{A\in\Xi}\chi_AP_{l,A}
\] the operator of of piecewise-polynomial approximation. Note that $rank(K_{\Xi,l})\leq l\cdot card(\Xi)\leq ln$.

\begin{thm}\label{thm:6}
With the above notation we have
\begin{equation}\label{eq:072}
\int_{\Omega}V(x)|u-K_{l,\Xi}u|^2dx \leq \kappa n^{-2l}\mathcal{J}(\Omega)\int_{\Omega}|\partial_r^lu|^2dx
\end{equation}
where $\kappa = \kappa (d,l,R_0,R)$
\end{thm}
{\bf Proof}: First
\begin{eqnarray}
\int_{\Omega}V(x)|u-K_{l,\Xi}u|^2dx = \sum_{A\in\Xi}\int_A V(x)|u-P_{l,A}u|^2dx \nonumber\\
\leq \sum_{A\in\Xi} \sup_{x\in A}\frac{|u(x)-P_{l,A}u(|x|)|^2}{\big[ 1+\log(R/|x|) \big]^{2i}} \cdot\int_{A}V(x)\big[ 1+\log(R/|x|) \big]^{2i}dx \nonumber\\
\leq \sum_{A\in\Xi} \sup_{x\in A}|u(x)-P_{l,A}u(|x|)|^2\cdot\mathcal{J}(A) \nonumber
\end{eqnarray}where we have used $\big[ 1+\log(R/|x|) \big]^{2i}\geq 1$ in the last inequality. So by equation (\ref{eq:071})
\begin{eqnarray}
\int_{\Omega}V(x)|u-K_{l,\Xi}u|^2dx \leq \kappa''\sum_{A\in\Xi}\int_A|\partial_r^lu|^2dx\cdot |A|^{2l-1}\cdot \mathcal{J}(A) \nonumber\\
\leq G(\mathcal{J},\Xi)\cdot \kappa''\cdot \sum_{A\in\Xi} \int_A|\partial_r^lu|^2dx \nonumber\\
\leq \kappa n^{-2l}\mathcal{J}(\Omega)\int_{\Omega}|\partial_r^lu|^2dx \nonumber
\end{eqnarray}where we used lemma (\ref{lemma5}) in the last line. $\blacksquare$
\\
\\
On the space $H_r = H_r(\Omega)$ endowed with the norm $(\int_{\Omega}|\nabla^lu|^2dx)^{1/2}$ we consider the quadratic form 
\[
b_V(u):= \int_{\Omega}V(x)|u(x)|^2dx
\]where $V(x)$ is radial and integrable when weighted with the weight $[1+\log(R/|x|)]^{2i}$. If $b_V$ is bounded on $H_r$, then it generates a bounded, self-adjoint, non-negative operator on $H_r$, which we will denote by $T_V$. By definition, for $f\in H_r$
\begin{eqnarray}
u=T_Vf &\iff&  u \in H_r; \nonumber \\
		&&\int_{\Omega}\nabla^lu\cdot\nabla^lwdx = \int_{\Omega}Vuwdx, \forall w\in H_r\nonumber
\end{eqnarray} As before we let
\[
n(\lambda;T_V)=N(-\lambda;-T_V) = \#\{\text{eigenvalues of $T_V$ that are $\geq \lambda$}\}
\]
\begin{thm}\label{thm:7}
Let $V(x)$ and $\Omega$ be as above. Then the operator $T_V$ is compact on $H_r$, and there exists a constant $C_4 =C_4(\Omega)$ such that for any $\lambda >0$
\[
n(\lambda;T_V)^{2l} \leq C_4\lambda^{-1} \int_{\Omega} V(x)\big[ 1+\log\big(R/|x|\big) \big]^{2i}dx
\]
\end{thm}
{\bf Proof}: We will actually prove an upper bound on for the eigenvalues of $T_V$ and use this to derive the estimate for $n(\lambda;T_V)$. Fix a natural number $n$. By Theorem \ref{thm:6} there exists a partition $\Xi$ of $\Omega$ into smaller annuli such that $card(\Xi)\leq n$ and for any radially symmetric $u\in H^l(\Omega)$ the estimate (\ref{eq:072}) holds. Let $\mathcal{F}:=ker(K_{\Xi,l})$. We have
\[
codim(\mathcal{F}) = rank(K_{\Xi,l}) \leq ln
\]
For $u\in \mathcal{F}$ we compute
\begin{eqnarray}
\langle T_Vu,u \rangle = b_V(u) &=& \int_{\Omega}V(x)|u|^2dx\nonumber\\
&=&\int_{\Omega}V(x)|u-K_{\Xi,l}u|^2dx \nonumber\\
&\leq& \kappa n^{-2l} \mathcal{J}(\Omega)\int_{\Omega}|\nabla^lu|^2dx\nonumber
\end{eqnarray}
Where Theorem (\ref{thm:6}) was applied in the last line. Thus by taking $n\to \infty$ we see that $T_V$ is the norm-limit of finite rank operators, hence compact. Also
\begin{equation}\label{eq:75}
\max_{u\in\mathcal{F}}\frac{\langle T_Vu,u \rangle_{H^l_0(\Omega)}}{\langle u,u \rangle_{H^l_0(\Omega)}} = \max_{u\in\mathcal{F}} \frac{\int_{\Omega}V(x)|u|^2dx}{\int_{\Omega}|\nabla^lu|^2dx} \leq \kappa n^{-2l} \cdot \mathcal{J}(\Omega) 
\end{equation}
So by Courant's Minimax Theorem for general symmetric compact operators we have
\[
\lambda_{nl+1} \leq \kappa n^{-2l}\cdot \mathcal{J}(\Omega)
\]
where $\lambda_j$ denotes the $j$-th of the positive operator $T_V$ in $H_r$. Hence there is a constant $C=C(\Omega)$ such that 
\begin{equation}\label{eq:76}
\lambda_n \leq C(\Omega) n^{-2l}\cdot \mathcal{J}(\Omega), \quad n=1,2, \ldots
\end{equation}
To prove the required estimate we proceed as follows
\begin{eqnarray}
n(\lambda;T_V) = \#\{n: \lambda_n \geq \lambda\}
&\leq& \#\{n: C(\Omega)\mathcal{J}(\Omega)n^{-2l}\geq \lambda \}\nonumber\\
&=& \#\{n: n\leq C(\Omega)^{1/2l}\mathcal{J}(\Omega)^{1/2l}\cdot \lambda^{-1/2l} \}\nonumber\\
&\leq& C(\Omega)^{1/2l}\mathcal{J}(\Omega)^{1/2l}\cdot \lambda^{-1/2l}\nonumber
\end{eqnarray}
raising both sides to the power $2l$ gives the result. $\blacksquare$\\
\\
{\bf Proof of Proposition \ref{prop:2}}: To finish the proof of Proposition \ref{prop:2} we again use the Birman-Swchinger Principle with $\mathcal{H} = L_{rad}^2(\Omega)$, $D(a) = \tilde{D}(a) = H_r$, $a(u)= \int_{\Omega}|\nabla^lu|^2dx$, $b(u) = \int_{\Omega}V|u|^2dx$, and $B = T_V$ to obtain
\[
N(0;(-\Delta)^l - V)^{2l} = n(1; T_V)^{2l} \leq C_{11} \mathcal{J}(\Omega).
\]
$\blacksquare$

\section{Problem (H), the radial problem with Hardy potential}\label{sec8}

Finally we consider problem (H). That is we seek weak radial solutions to
\[
\tag{H}
\begin{cases}
(-\Delta)^lu + \frac{\bar{b}}{|x|^{2l}}u= g(x,u) + \varphi(x) \quad \mbox{in $B_R(0)$} 
\\
\\
\bigg(\frac{\partial}{\partial \nu}\bigg)^ju\Bigg|_{\partial B_R(0)} = 0, \quad j = 0, \dots, l-1
\end{cases}
\]where $\bar{b}>0$, and $g$ and $\varphi$ are radially symmetric in their dependence on $x$. Notice that the Hardy exponent is the critical one. That is the corresponding Hardy inequality in $H_0^l(B_R)$ doesn't hold in general. Define the space $\mathcal{C}_r$ as the set of $u\in C^\infty_0(B_R\setminus\{0\})$ with $u(x)=u(|x|)$ a.e. Let
\[
\mathcal{H}_r := \text{the closure of $\mathcal{C}_r$ in the following norm}
\]

\begin{equation}\label{eq:77}
||u||_{\mathcal{H}_r}^2 := \int_{B_R}|\nabla^lu|^2dx + \bar{b}\int_{B_R}\frac{|u|^2}{|x|^{2l}}dx
\end{equation}Note that clearly $\mathcal{H}_r \hookrightarrow H^l_0(B_R)$ continuously, something to keep in mind. We will prove the existence of solutions by looking for critical points of the functional
\begin{equation}
I_1(u) := \frac{1}{2}\int_{B_R}|\nabla^lu|^2dx +\frac{\bar{b}}{2}\int_{B_R}\frac{|u|^2}{|x|^{2l}}dx- \int_{B_R} G(x,u)dx - \int_{B_R}\varphi u dx\nonumber
\end{equation} on $\mathcal{H}_r$, where $G(x,u)=\int_0^ug(x,t)dt$. Because $g$ and $\varphi$ are radially symmetric in their explicit dependence on $x$, the critical points of this functional on $\mathcal{H}_r$ are also critical points on $H^l_0(B_R)$. Hence they are generalized solutions for (H). Since $G(x,u)$ has the same growth restrictions in this problem as it did in (P) all conditions in Theorem \ref{thm:3} and in Tanaka's theorem follow exactly the same reasoning as they did earlier. Since $\bar{b}>0$ Palais-Smale sequences are bounded in the norm of $\mathcal{H}_r$, instead of merely in the norm of $H_0^l$. In any case the argument proceeds as before and we need to contradict the possibility that the minimax values $c_n$ of the unperturbed functional satisfy equation (\ref{eqApp29})
\begin{equation}
c_n \leq An\big[\ln(n)\big]^{1/\alpha_1} \quad \text{for $n>n_0$}.\nonumber
\end{equation}
Earlier in section \ref{sec4} we define the smoother functional $J(u)$ by equation (\ref{eq:31}), and its minimax levels
\begin{equation}
b_n := \inf_{g\in\Gamma_n}\sup_{u\in g(E_n)}J(u).\nonumber
\end{equation}
We do exactly the same thing here, except $J(u)$ is now defined as
\[
J(u) = \frac{1}{2}\int_{B_R}|\nabla^lu|^2dx +\frac{\bar{b}}{2}\int_{B_R}\frac{|u|^2}{|x|^{2l}}dx - \int_{B_R}H(u)dx.
\]
That is, the form $\frac{1}{2}\int_{B_R}|\nabla^lu|^2dx$ is replaced by the from $\frac{1}{2}\int_{B_R}|\nabla^lu|^2dx +\frac{\bar{b}}{2}\int_{B_R}\frac{|u|^2}{|x|^{2l}}dx$. Since $J(u)\leq I_0(u)$ by construction, we have $b_n \leq c_n$. So it will suffice to obtain a good lower bound on $b_n$. By Tanaka \cite{Tanaka89} Theorem B, there exists a sequence $u_n$ such that

\begin{enumerate}
\item[i)] $J(u_n)\leq b_n$

\item[ii)] $J'(u_n) = 0$

\item[iii)] $n\leq index_0J''(u_n)$\\
\end{enumerate} Therefore 
\[
n \leq index_0J''(u) = N_{-}((-\Delta)^l + b|x|^{-2l} - H''(u_n(x)))
\]The eigenvalue estimate we apply is the following result from \cite{Laptev-Netrusov}
\begin{lemma}[Laptev-Netrusov]
Consider the unbounded linear operator $H:=(-\Delta)^l + b|x|^{-2l} - V(x)$ acting on $L^2(B_R)$, $B_R$ the ball of radius $R$ in $\mathbb{R}^{2l}$. Suppose $V(x)=V(|x|)\geq 0$ and $V(x)\in L^1(B_R)$. Then
\[
N_{-}((-\Delta)^l + b|x|^{-2l} - V(x)) \leq C(l,b)\int_{B_R} V(x)dx
\]
\end{lemma} So
\[
n \leq C\int_{B_R} H''(u_n(x))dx \leq C\int_{B_R}(|u_n|+1)^{2\alpha_2-2}e^{(u_n^2+1)^b}dx
\]Since $u_n$ is a critical point of $J$ the last term in the above line can be controlled by $J(u_n)$. Like in section \ref{sec4}, this gives
\[
Cn [\log(n)]^{\gamma} \leq J(u_n)\leq b_n \leq c_n
\] for sufficiently large $n$, where this time $\gamma = 2/\alpha_2-1$. Thus if $1/\alpha_1 < 2/\alpha_2-1$ this contradicts the alleged upper bound, and proves that $I_1$ an unbounded sequence of radial critical points.
\\
\\

{\small {\bf Acknowledgements} The author would like to thank Marcello Lucia, Mythily Ramaswamy, and Kyril Tintarev for some helpful discussions.}


\begin{thebibliography}{1}

\bibitem{DRAdams88} D. R. Adams, \emph{A sharp inequality of J. Moser for higher order derivatives}, Ann. Math. 128, (1988), 385-398.

\bibitem{Adams75} R. A. Adams, \emph{Sobolev Spaces}, Academic Press, New York-San Francisco-London, 1975.


\bibitem{AmbrosettiRabinowitz1} A. Ambrosetti, P. H. Rabinowitz, \emph{Dual variational methods in critical point theory and applications}, L. Funct. Anal. 14 (1973), 349-381.

\bibitem{Ba-Be1} A. Bahri, H. Berestycki, \emph{A perturbation method in critical point theory and applications}, TAMS 267, (1981), 1-32.

\bibitem{Ba-Be2} A. Bahri, H. Berestycki, \emph{Forced vibrations of superquadratic Hamiltonian systems}, Acta Math. 152, (1984), 143-197.

\bibitem{Ba-Lions} A. Bahri, P. L. Lions, \emph{Morse index of some min-max critical points I. Applications to multiplicity results}, Comm. Pure Appl. Math. 41, (1988), 1027-1037.

\bibitem{BirmanSolomyak2} M. Sh. Birman, M. Z. Solomyak, \emph{Quantitative analysis in Sobolev imbedding theorems and applications to spectral theory}, Tenth Math. Summer School (Katsiveli/Nalchik, 1972), Izdanie Inst. Mat. Akad. Nauk. Ukrain. SSR, Kiev, 1974, pp. 5-189; English transl. Amer. Math. Soc. Transl. (2) 114 (1980).

\bibitem{BirmanSolomyakBook} M. Sh. Birman, M. Z. Solomyak, \emph{Spectral Theory of Self-Adjoint Operators in Hilbert Space}, Mathematics and its Applications 5, Reidel Publ. Company, Dordrecht (1987).

\bibitem{BirmanSolomyak3} M. Sh. Birman, M. Z. Solomyak, \emph{Estimates for the number of negative eigenvalues of the Schrodinger operator and its generalizations}, Adv. in Sov. Math., Vol. 7 (M. Birman, ed.), AMS, 1991, 1-55.

\bibitem{Bolle99} P. Bolle, \emph{On the Bolza problem}, J. Differential Equations, 152, (1999), 274-288.

\bibitem{Bo-Gh-Te} P. Bolle, N. Ghoussoub, H. Tehrani, \emph{The multiplicity of solutions in non-homogeneous boundary value problems}, Manuscripta Math. 101, (2000), 325-350.

\bibitem{Clapp-Ding-Hernandez} M. Clapp, Y. Ding, S. Hern\`andez-Linares, \emph{Strongly indefinite functionals with perturbed symmetries and multiple solutions of non symmetric elliptic systems}, Electron. J. Differential Equations 100, (2004), 1-18.

\bibitem{Cwikel77} M. Cwikel, \emph{Weak type estimates for singular values and the number of bound states of Schrodinger operators}, Ann. Math. 106, (1977), 93-100.

\bibitem{Gazzola-Grunau-Sweers} F. Gazzola, H.-Ch. Grunau, G. Sweers, \emph{Polyharmonic boundary value problems}, Springer-Verlag, Berlin Heidelberg, 2010.

\bibitem{Kras-Rutic} M. A. Krasnosel'skii and Ya. B. Rutickii, \emph{Convex functions and Orlics spaces}, Moscow, Fizmatgiz, 1958(Russian); English transl.: P. Nordhoff, Groningen, 1961.

\bibitem{La-Mu-Sq} S. Lancelotti, A. Musesti, M. Squassina, \emph{Infinitely many solutions for polyharmonic elliptic problems with broken symmetries}, Math. Nachr. 253, (2003), 35-44.

\bibitem{Laptev-Netrusov} A. Laptev, Yu. Netrusov, \emph{On the negative eigenvalues of a class of Schrodinger operators}, Amer. Math. Soc. Transl. 189 (1999), 173 - 186.

\bibitem{Lieb76} E. Lieb, \emph{The number of bound states of one-body Schrodinger operators and the Weyl problem}, Bull. Amer. Math. Soc. 82, (1976), 751-753. 

\bibitem{Mitidieri} E. Mitidieri, \emph{A simple approach to Hardy's inequalities}, Math. Notes 67, (2000), 479-486.

\bibitem{Moser71} J. Moser, \emph{A sharp form of an inequality by N. Trudinger}, Ind. Univ. Math. J. 20, (1971), 1077-1092.

\bibitem{Ni1} W.-M. Ni, \emph{A nonlinear Dirichlet problem on the unit ball and its applications}, Indiana Univ. Math. J. 31, (1982), 801-807.

\bibitem{Palais79} R. Palais. \emph{The principle of symmetric criticality}, Comm. Math. Phys. 69 (1979), no. 1, 19--30.

\bibitem{Pohozaev65} S. I. Pohozaev, \emph{The Sobolev embedding in the case $pl=n$}, Proceedings of the Technical Scientific Conference on Advances of Scientific Research 1964-1965. Mathematics Section, Moskov. Energet. Inst., Moscow (1965), 158-170.

\bibitem{Rabinowitz82} P. H. Rabinowitz, \emph{Multiple critical points of perturbed symmetric functionals}, Trans. Amer. Math. Soc. 272, (1982), 753-769.

\bibitem{Rabinowitz86} P. H. Rabinowitz, \emph{Minimax Methods in Critical Point Theory with Applications to Differential Equations}, CBMS, AMS 65, (1986).

\bibitem{RaoRen} M. M. Rao, Z. D. Ren, \emph{Theory of Orlicz spaces}, Marcel Dekker, New York, 1991.

\bibitem{Rozen72} G. V. Rozenblum, \emph{The distribution of the discrete spectrum for singular differential operators}, Dokl. Akad. Nauk SSSR 202 (1972), 1012-1015; English transl. in: Soviet Math. Dokl. 13 (1972).

\bibitem{Rozen76}G. V. Rozenblum, \emph{Distribution of the discrete spectrum for singular differential operators}, Izv. Vyssh. Uchebn. Zaved. Mat. No. 1 (1976), 75-86: English transl. in: Soviet Math. (Iz. VUZ.) 20, (1976).

\bibitem{Solomyak94} M. Z. Solomyak, \emph{Piecewise-polynomial approximation of functions from $H^l((0,1)^d)$, $2l=d$, and applications to the spectral theory of the Schrodinger Operator}, Isr. J. Math. 86 (1994), 253-275



\bibitem{Stuwe80} M. Struwe, \emph{Infinitely many critical points for functionals which are not even and applications to superlinear boundary value problems}, Manuscripta Math. 32, (1980), 335-364.

\bibitem{Struwe08} M. Struwe, \emph{Variational Methods. Applications to Nonlinear Partial Differential Equations and Hamiltonian Systems}, Fourth Edition, Springer-Verlag, Heidelberg, 2008.

\bibitem{Sugimura94} K. Sugimura,\emph{Existence of infinitely many solutions for a perturbed elliptic equation with exponential growth}, Nonlin. Anal., Theory, Meth. Appl., 22 (1994), 277-293.

\bibitem{Tanaka89} K. Tanaka, \emph{Morse indices at critical points related to the symmetric mountain pass theorem and applications}, Comm. Partial Differential Equations 14, (1989), 99-128.

\bibitem{Tarsi08} C. Tarsi, \emph{Perturbation from symmetry and multiplicity of solutions for elliptic problems with subcritical exponential growth in $\mathbb{R}^2$}, Comm. Pure Appl. Anal. 7 (2), (2008), 445. 

\bibitem{Trudinger67} N. Trudinger, \emph{On imbeddings into Orlicz spaces and some applications}, J. Math. Mech., 17, (1967), 473-484.

\bibitem{Weidmann} J. Weidmann, \emph{Linear operators in Hilbert spaces}, First edition, Springer-Verlag, New York, 1980.



\end{thebibliography}
\end{document}